\documentclass[12pt]{amsart}
\let\cal\mathcal
\def\A{\mathcal A}
\def\B{\mathcal B}
\def\C{\mathcal C}
\def\D{\mathcal D}
\def\I{\mathcal I}
\def\K{\mathcal K}

\def\M{\mathcal M}
\def\X{\mathcal X}
\def\P{\mathcal P}
\def\Y{\mathcal Y}
\def\H{\mathcal H}
\def\U{\mathcal U} 
\def\ZZ{\mathbb Z} 
\def\NN{\mathbb N} 
\def\CC{\mathbb C} 
\def\RR{\mathbb R} 
\def\DD{\mathbb D} 
\def\ourimath{i}
\def\ourjmath{j}
\def\id{\rm id}
\def\iso{\cong}
\def\ptp{\hat{\otimes}}
\newcommand{\Tot}{\mathop{\rm Tot}}
 
\newcommand{\frechet}{Fr\'{e}chet}
\newcommand{\fr}{{\mathcal F} r}
\newcommand{\im}{\mathop{\rm Im}}
\newcommand{\larray}{\left(\begin{array}{cc}}
\newcommand{\rarray}{\end{array}\right)}

\newtheorem{theorem}{\sc Theorem}[section]
\newtheorem{lemma}[theorem]{\sc Lemma}
\newtheorem{proposition}[theorem]{\sc Proposition}
\newtheorem{corollary}[theorem]{\sc Corollary}
\newtheorem{remark}[theorem]{\sc Remark}
\newtheorem{example}[theorem]{\sc Example}
\newtheorem{definition}[theorem]{\sc Definition}
\begin{document}

\title{The K\"unneth formula for 
nuclear $DF$-spaces and Hochschild cohomology}

\author{Zinaida A. Lykova}

\address{School of Mathematics and Statistics, University of Newcastle,\\
 Newcastle upon Tyne, NE1 7RU, UK~ {\rm (Z.A.Lykova@newcastle.ac.uk)}}

\date{September 2007}

\begin{abstract} We consider complexes $(\X, d)$ of nuclear Fr\'echet spaces 
and continuous boundary maps $d_n$ with closed ranges and prove that,
up to topological isomorphism, $\;(H_{n}(\X, d))^*$ $\iso$ $H^{n}(\X^*,d^*),$ where 
$(H_{n}(\X, d))^*$ is the strong dual space of the homology group of 
$(\X, d) $ and $\; H^{n}(\X^*,d^*)$ is the cohomology group of the strong 
dual complex $ (\X^*,d^*)$. We use this result to  establish  the existence 
of topological isomorphisms in the  K\"{u}nneth formula  for the 
cohomology  of complete nuclear $DF$-complexes  and in the  
K\"{u}nneth formula  for continuous Hochschild cohomology
of  nuclear  $\hat{\otimes}$-algebras which are Fr\'echet spaces or 
 $DF$-spaces for which all boundary maps of the standard homology complexes have closed ranges.  
We describe explicitly continuous Hochschild and cyclic  
cohomology groups of certain tensor products of
$\hat{\otimes}$-algebras which are Fr\'echet spaces or 
 nuclear  $DF$-spaces.\\

\end{abstract}


\thanks{I am indebted to the Isaac Newton Institute for Mathematical 
Sciences at Cambridge for hospitality and for generous financial 
support from the programme on Noncommutative Geometry while 
this work was carried out.}


\maketitle
	
\markboth{Z.~A.~Lykova}{The K\"unneth formula for nuclear $DF$-spaces}

\section{Introduction}
K\"unneth formulae for bounded chain  complexes  $\X$ and $\Y$ 
of \frechet\  and Banach spaces and continuous boundary maps  
with closed ranges  were established,  under 
certain topological assumptions, in \cite{Ka,GLW,GLW2}. Recall that in the category of nuclear \frechet\ spaces short exact sequences are topologically 
pure and objects are strictly flat, and so the K\"unneth formula 
can be used for calculation of continuous Hochschild homology 
$\H_n(\A \hat{\otimes} \B, X \hat{\otimes} Y)$ if 
boundary maps of the standard homology complexes have closed ranges. 
To compute the continuous cyclic-type Hochschild cohomology  
of \frechet\ algebras one has to deal with complexes of complete 
$DF$-spaces.  In the recent paper \cite{Ly5} the author showed 
that,  for a continuous morphism $ \varphi: \X
\rightarrow \Y$ of complexes of complete nuclear $DF$-spaces, a surjective map of cohomology groups
$H^n(\varphi):  H^n(\X) \rightarrow H^n(\Y)$ is automatically open.
In this paper we establish relations between topological properties of the  homology of complexes of Fr\'echet spaces and of the cohomology of their strong dual complexes. We use these properties  to show the existence of a topological isomorphism in the  K\"{u}nneth 
formula for complexes of complete nuclear $DF$-spaces and continuous 
boundary maps  with closed ranges and thereby
 to describe explicitly the  continuous Hochschild and 
cyclic  homology and cohomology of $\A \hat{\otimes} \B$ for 
certain  $\hat{\otimes}$-algebras $\A$ and $\B$ which 
are Fr\'echet spaces or  nuclear  $DF$-spaces.

 In Theorem \ref{H(dualX)=dualH(X)} and 
Corollary \ref{H(dualX) =dualH(X)-DF}, for a complex of nuclear 
Fr\'echet spaces or of complete nuclear $DF$-spaces $(\X, d) $ and continuous boundary maps $d_n$ with closed ranges, 
we establish that there is a topological isomorphism,
$ (H_{n}(\X, d))^* \iso H^{n}(\X^*,d^*),$ where 
$ (H_{n}(\X, d))^*$ is the dual space of the homology group of 
$(\X, d) $ and $ H^{n}(\X^*,d^*)$ is the cohomology group of the dual
complex $ (\X^*,d^*)$. 

In  Theorem \ref{NuclearCase-DF} and Theorem \ref{NuclearCase}, for  
bounded chain complexes $\X$ and $\Y$ of complete nuclear $DF$-spaces 
or of nuclear \frechet\ spaces such that all boundary maps 
have closed ranges, we prove that, up to  topological  isomorphism, 
 $$ H_n(\X \ptp \Y) \iso\displaystyle{\bigoplus_{m+q=n}} 
[H_m(\X) \ptp H_q(\Y)] \;\;{\rm and}$$ 
$$ H^n((\X \ptp \Y)^*) \iso
\displaystyle{\bigoplus_{m+q=n}} [H^m(\X^*) \ptp H^q(\Y^*)].$$
In Corollary \ref{Kunneth-Ban-Fr}, for  
bounded chain complexes  $(\X, d_{\X}) $ of Banach spaces and $(\Y, d_{\Y}) $ of \frechet\ spaces 
such that all boundary maps have closed ranges, and $H_n(\X)$ and 
${\rm Ker}~ (d_{\X})_n$ are strictly flat in ${\mathcal Ban}$
for all $n$, we prove that, 
up to  topological  isomorphism, 
 $$ H_n(\X \ptp \Y) \iso\displaystyle{\bigoplus_{m+q=n}} 
[H_m(\X) \ptp H_q(\Y)]$$
and, up to isomorphism of linear spaces,
$$ H^n((\X \ptp \Y)^*) \iso \displaystyle{\bigoplus_{m+q=n}} [H_m(\X) \ptp H_q(\Y)]^*.$$

The K\"unneth formulae for the continuous Hochschild homology  
$\H_n(\A \hat{\otimes} \B, X \hat{\otimes} Y)$ and cohomology 
$H^n(({\mathcal C}_{\sim} (\A \ptp \B, X \ptp Y))^*)$ are proved in
 Theorem \ref{Hochschild-(co)hom} for the underlying category of 
complete nuclear $DF$-spaces and for the underlying category of nuclear \frechet\ spaces.
In these  underlying categories, for unital  $\hat{\otimes}$-algebras $\A$ and $\B$, for  a unital $\A$-$\hat{\otimes}$-bimodule $X$ and a unital $\B$-$\hat{\otimes}$-bimodule $Y$, under the assumption that all boundary maps of the standard homology complexes ${\mathcal C}_{\sim}(\A, X)$
and ${\mathcal C}_{\sim}(\B, Y)$ have closed ranges, we show that, up 
to topological isomorphism, 
$$ \H_n(\A \ptp \B, X \ptp Y) \iso
\displaystyle{\bigoplus_{m+q=n}} [\H_m(\A , X) \ptp 
\H_q(\B, Y)] 
\;\;{\rm and}$$
$$H^n(({\mathcal C}_{\sim} (\A \ptp \B, X \ptp Y))^*) \iso $$
$$ \displaystyle{\bigoplus_{m+q=n}}
[H^m(({\mathcal C}_{\sim} (\A, X))^*) 
\ptp H^q(({\mathcal C}_{\sim} (\B, Y))^*)]. $$
In Theorem \ref{Hochschild-(co)hom-Ban-Fr} we prove the K\"unneth formulae for 
the continuous Hochschild homology groups of Banach and \frechet\ algebras under some topological assumptions.
In Section 6 we describe explicitly the continuous cyclic-type homology
and cohomology of certain tensor products of $\hat{\otimes}$-algebras
which are Banach or  \frechet\  or nuclear  \frechet\  or nuclear $DF$-spaces.

\section{Definitions and notation}

 We recall some notation and terminology used in homology 
and in the theory of topological algebras. Homological theory can be found in any relevant textbook, for instance, 
MacLane \cite{Mac}, Loday \cite{Loday} for the pure algebraic case 
and  Helemskii \cite{He0} for the continuous case.

Throughout the paper $\hat{\otimes}$ is the projective tensor product of 
complete locally convex spaces. By $ X^{\ptp n}$ we mean the $n$-fold projective tensor power $ X\ptp \dots \ptp X $ of $X$ and ${\rm id}$
denotes the identity operator. 

We use the notation ${\mathcal Ban}$,  $\fr$  and $\mathcal{LCS}$
for the categories whose objects are Banach spaces, Fr\'{e}chet 
spaces and complete Hausdorff locally convex
spaces respectively, and whose morphisms in all cases are continuous
linear operators. 
For topological homology theory it is important to 
find a suitable category for the underlying spaces of the algebras 
and modules. In \cite{He0} Helemskii constructed homology theory for 
the following categories $\Phi$ of underlying spaces, for which he used the
notation $(\Phi, \hat{\otimes})$.

\begin{definition}\label{category}(\cite[Section II.5]{He0}) A {\em suitable category} for underlying spaces of the algebras and modules is
an arbitrary complete subcategory 
$\Phi$ of $\mathcal{LCS}$ having 
the following properties:

(i) if  $\Phi$ contains a space, it also contains all
those spaces topologically isomorphic to it;

(ii) if  $\Phi$ contains a space, it also contains 
any of its closed subspaces and the completion of any its Hausdorff
quotient spaces;

(iii)  $\Phi$ contains the direct sum and the
projective tensor product of any pair of its spaces;

(iv) $\Phi$ contains  ${\mathbb C}$. 
\end{definition}

Besides ${\mathcal Ban}$,  $\fr$  and $\mathcal{LCS}$ important
examples of suitable categories $\Phi$ are the categories of
complete nuclear spaces \cite[Proposition 50.1]{Tr}, 
 nuclear \frechet\ spaces and complete nuclear $DF$-spaces \cite{Ly5}.

By definition a $\hat{\otimes}$-algebra is a complete Hausdorff  
locally convex  algebra with jointly continuous multiplication. 
A left  $\hat{\otimes}$-module $X$ over a $\hat{\otimes}$-algebra $\A$ 
is a complete Hausdorff  locally convex space  $X$ together with the 
structure of a left $\A$-module such that the map 
$\A \times X \to X$, $ (a,x) \mapsto a \cdot x $ is jointly continuous. 
For a  $\hat{\otimes}$-algebra 
$\A$, $\hat{\otimes}_\A$ is the projective tensor product over $\A$
of   left and right  $\A$-$\hat{\otimes}$-modules (see \cite{GrA},
\cite[II.4.1]{He0}). The category of  left  [unital] 
$\A$-$\hat{\otimes}$-modules is denoted by 
$\A$-{\rm mod} [$\A$-{\rm unmod}] and the category 
of  [unital]  $\A$-$\hat{\otimes}$-bimodules 
is denoted by  $\A$-{\rm mod}-$\A$ [$\A$-{\rm unmod}-$\A$]. 

 Let ${\K}$ be a category.  A {\it  chain complex} ${\cal X}_{\sim} $ in the category  ${\K}$ is a  sequence of $ X_n \in {\K}$ and  morphisms 
$d_n $ (called {\em boundary maps})
\[
 \dots \leftarrow X_n \stackrel {d_n} {\leftarrow}  X_{n+1} 
\stackrel {d_{n+1}} {\leftarrow} X_{n+2}
 \leftarrow  \dots
\]
such that $d_n \circ d_{n+1} = 0 $ for every $n$.
 The {\it cycles} are the elements of 
$$Z_n(\X) = {\rm Ker}~( d_{n-1}:
X_n \rightarrow  X_{n-1}).$$ 
The   {\it boundaries} are the elements of 
$$B_n(\X) = {\rm Im}~( d_n: X_{n +1} \rightarrow  X_n).$$ 
The relation $d_{n-1} \circ d_n = 0 $
implies $B_n(\X) \subset Z_n(\X).$ 
The {\it homology groups} are defined by $$H_n({\cal X}_{\sim}) = 
Z_n(\X)/B_n(\X).$$
As usual, we will often drop the subscript $n$ of $d_n$. If there 
is a need to distinguish between various boundary maps on various 
chain complexes, we will use subscripts, that is, we will denote 
the boundary maps on $\X$ by $d_{\X}$.
 A chain complex $\X$ is called {\it bounded} if $X_{n} = \{ 0 \}$ 
whenever $n$ is less than a certain fixed integer $N \in \ZZ$.

Given  $E \in \K$ and a chain complex $(\X, d)$ in $\K$, we can 
form the chain complex $E \ptp \X$ of the locally convex spaces 
$E\ptp X_n$ and boundary maps $\id_E \otimes d$.  
Definitions of the totalization $\Tot(\overline{\M})$ of a bounded 
bicomplex $\overline{\M}$ and the tensor product $\X\ptp\Y$ of bounded 
complexes $\X$ and $\Y$ in $\fr$
  can be found in \cite[Definitions II.5.23-25]{He0}. Recall that  
$\X\ptp\Y \stackrel{\rm def}{=} {\rm Tot} (\overline{\X\ptp\Y})$ of  a 
bounded bicomplex  $\overline{\X\ptp\Y}$.

We recall here the definition of a strictly flat 
locally convex  space in a suitable category $\Phi$ 
which is equivalent to that given in \cite[Chapter 
VII]{He0}. Note that it can be seen as a special case of the 
corresponding notion for $\hat{\otimes}$-modules, where the 
$\hat{\otimes}$-algebra is taken to be the complex numbers $\CC$.

\begin{definition}\label{strict_flat} A locally convex space $G \in \Phi$
is {\em strictly flat} in $\Phi$ if
for every short exact sequence 
$$ 0 \to X \to Y \to Z \to 0$$
of locally convex spaces from $\Phi$ and continuous
linear operators, the short sequence
$$ 0 \to G\hat\otimes X \to G\hat\otimes Y \to G\hat\otimes Z \to 0$$
is also exact.
\end{definition}

\begin{example} {\rm (i)} Nuclear Fr\'echet spaces are strictly flat 
in $\fr $ \cite[Theorems A.1.5 and A.1.6]{EsPu}. 
{\rm (ii)} Finite-dimensional Banach spaces and $L^1(\Omega, \mu)$ are 
strictly flat in ${\mathcal B}an$ \cite[Theorem III.B.2]{Woj} and 
in $\fr$ \cite[Proposition 4.4]{Pir3}.
\end{example}

If $E$ is  a topological vector space $E^*$ denotes its dual space of
continuous linear functionals.
For a subset $V$ of $E$, the {\it polar} of  $V$ is
\[V^0= \{ g  \in E^*: |g(x)| \le 1 \; {\rm for\; all}\; x \in V \}.\]
 Throughout the paper, $E^*$ will always be equipped with 
the strong topology unless otherwise stated. The  {\it strong topology}
is defined on $E^*$ by taking as a basis of 
neighbourhoods of $0$ the family of polars of all bounded subsets of $E$; see \cite[II.19.2]{Tr}.

Let $\A$ be a $\hat{\otimes}$-algebra.
 A complex of  $\A$-$\hat{\otimes}$-modules and their morphisms
is called {\it  admissible} if it 
splits as a complex in  $\mathcal{LCS}$ \cite[III.1.11]{He0}.
 A complex of  $\A$-$\hat{\otimes}$-modules and their morphisms
is called {\it weakly admissible} if its strong dual complex  
splits.

 For $Y \in \A$-{\rm mod}-$\A$ a complex 

\vspace*{0.2cm}
\hspace{1.5cm}
$0  \longleftarrow Y
\stackrel {\varepsilon} { \longleftarrow}  P_0 
 \stackrel {\phi_0} { \longleftarrow} P_1 \stackrel {\phi_1} 
{ \longleftarrow} P_2 \longleftarrow  \cdots  
\hfill {(0 \leftarrow Y \leftarrow \P)} $ 
\vspace*{0.2cm}
\newline is called a {\it   projective
resolution  of  $Y$ in  $\A$-{\rm mod}-$\A$ } if it is admissible
 and all the modules in ${\mathcal P}$ are  
projective in $\A$-{\rm mod}-$\A$ \cite[Definition III.2.1]{He0}. 

 For any $\hat{\otimes}$-algebra $\A$, not necessarily unital, 
$\A_+$ is the  $\hat{\otimes}$-algebra obtained by 
adjoining an identity to $\A$.
 For a $\hat{\otimes}$-algebra $\A$, the algebra 
$\A^e = \A_+ \hat{\otimes} \A_+^{op}$
is called the {\it 
 enveloping algebra of $\A$}, where $\A_+^{op}$ is the 
{\it  opposite algebra}
of $\A_+$ with  multiplication $a \cdot b = ba$. 
 
 A module $Y \in \A$-{\rm mod} is called 
{\it flat } if for any admissible complex ${\X}$ of right 
 $\A$-$\hat{\otimes}$-modules the complex ${\X} \hat{\otimes}_{\A} Y$
is exact. A module $Y \in  \A$-{\rm mod}-$\A$ is called 
{\it flat } if for any admissible complex ${\X}$ of  
 $\A$-$\hat{\otimes}$-bimodules the complex 
${\X} \hat{\otimes}_{\A^e} Y$ is exact.

 For $Y, X  \in \A$-{\rm mod}-$\A$, we shall denote 
by  $ {\rm Tor}_{n}^{\A^e}(X, Y)$ the $n$th homology 
of the complex  $ X \ptp_{\A^e}{\mathcal P}$, where 
$0 \leftarrow Y \leftarrow {\mathcal P}$ is a  projective resolution 
of $Y$ in $\A$-{\rm mod}-$\A$,
  \cite[Definition III.4.23]{He0}.

\begin{definition}\label{topol_pure}~
 A short exact sequence of locally convex spaces from  
$\Phi$ and continuous operators
$$ 0 \rightarrow Y \stackrel {i} {\rightarrow}    Z  \stackrel {j} 
{\rightarrow} W  \rightarrow  0$$
is called {\em topologically pure in $\Phi$} if for every $X \in
\Phi$ the sequence
$$ 0 \rightarrow  X \hat\otimes Y \stackrel {{\rm id}_X \hat{\otimes} i} 
{\rightarrow}    X \hat\otimes Z  
 \stackrel {{\rm id}_X \hat{\otimes} j} {\rightarrow}    X \hat\otimes W
 \rightarrow  0 $$
is exact. 
\end{definition}

By \cite[II.1.8f and Remark after II.1.9]{CLM}, an  extension 
of Banach spaces is topologically pure in ${\mathcal B}an$
if and only if it is weakly admissible in ${\mathcal B}an$.
In the category of Fr\'{e}chet  spaces the situation 
with topologically pure extensions is more interesting.
Firstly, it is known that
extensions of nuclear Fr\'{e}chet  spaces are topologically pure
(see \cite[Theorems A.1.6 and A.1.5]{EsPu}. 
Note that nuclear Fr\'{e}chet   spaces are reflexive, and therefore 
a short sequence of nuclear Fr\'{e}chet   
spaces is weakly admissible if and only if it is admissible.
It is shown in \cite[Lemma 2.4]{Ly4}  that in $\fr$ the weak 
 admissibility of an extension
implies topological purity of the extension, but 
is not equivalent to the topological purity of the extension 
 \cite[Section 2]{Ly4}. Recall that  extensions of Fr\'{e}chet 
 algebras $\; 0 \to Y \to Z \to W \to 0\;$
such that $ Y$ has a left or right bounded approximate identity  
are topologically pure \cite[Lemma 2.5]{Ly4}. 

\section{Topological isomorphism between $(H_n(\X, d))^*$ and $H^n(\X^*, d^*)$ 
in the category of complete nuclear $DF$-spaces}

$DF$-spaces were introduced by A. Grothendieck in \cite{GrA1}. It is well known that the strong dual of a Fr\'{e}chet space is a complete $DF$-space and that nuclear Fr\'{e}chet spaces and complete nuclear $DF$-spaces are reflexive \cite[Theorem 4.4.12]{Pi}. Moreover, the correspondence 
$E \leftrightarrow E^*$ establishes a one-to-one relation between the nuclear 
Fr\'{e}chet spaces and complete nuclear $DF$-spaces 
\cite[Theorem 4.4.13]{Pi}. It is known that there exist closed linear subspaces of $DF$-spaces that are not $DF$-spaces. For nuclear spaces, however, we have the following.

\begin{lemma}\label{closed_suspace_DFspaces} \cite[Proposition 5.1.7]{Pi}~
 Each closed linear subspace $F$ of the strong dual of
a nuclear Fr\'echet space $E$ is also the strong dual of a nuclear 
Fr\'echet space.
\end{lemma}

In a locally convex space a subset is called a {\it barrel} if it 
is absolutely convex, absorbent and closed. Every locally
convex space has a neighbourhood base consisting of barrels. 
A locally convex space is called  a {\it barrelled} space
or a {\it $t$-space} if every barrel is a neighbourhood \cite{RR}. 
A Hausdorff barrelled locally convex space with the further property 
that its closed bounded subsets are compact is called a Montel space.
In particular, each nuclear Fr\'echet space is an $FM$-space, that is, a 
Fr\'echet space that is Montel.

Let $E$ and $F$ be locally convex topological vector spaces and let 
$T:E \rightarrow F$ be  a continuous linear
 operator. If $T$ is open, then it is also called a {\it homomorphism}.
Let $E^*$ and $F^*$ have the strong dual topologies. If 
$T^*:F^* \rightarrow E^*$ is an  open continuous linear
 operator then $T^*$ is called a {\it  strong homomorphism}.

\begin{lemma}\label{KerT*barrelled} 
Let $X$ and $Y$ be nuclear Fr\'echet spaces  and let $T:X \rightarrow Y$ be
a continuous linear operator such that ${\rm Im}~ T$ is closed. Then
 up to topological isomorphism, 
${\rm Im}~ T^* \cong (X/{\rm Ker}~ T)^*$ is 
the strong dual of a nuclear Fr\'echet space and hence is barrelled.
\end{lemma}

{\it Proof.} ~  In view of the algebraic identification
${\rm Im}~ T^* \cong (X/{\rm Ker}~ T)^*$, the adjoint of the quotient map
$X  \rightarrow X/{\rm Ker}~ T$ is the inclusion map ${\rm Im}~ T^*
\rightarrow X^*$. Since  nuclear Fr\'echet spaces are in particular 
$FM$-spaces  and the  adjoint of a homomorphism 
between $FM$-spaces is  a strong homomorphism, the relative topology and
the strong dual topology of ${\rm Im}~ T^*$ (as the dual space of 
$ X/{\rm Ker}~ T$) coincide. Hence, up to topological isomorphism, 
${\rm Im}~ T^* \cong (X/{\rm Ker}~ T)^*$.
By \cite[Corollary IV.3.1]{RR}, the strong dual of a 
nuclear Fr\'echet space is barrelled.
\hfill\mbox{$\Box$}\vspace{2mm}

Further we will need the following version of the open mapping theorem.

\begin{corollary}\label{open-map-strong-dual-FM}
\cite[Corollary 3.6]{Ly5}  Let $E$ and $F$
be  nuclear Fr\'echet spaces and let $E^*$ and $F^*$ be the strong duals of  $E$ and $F$ respectively.
Then a  continuous linear operator  $T$ of $E^*$ onto $F^*$ is open.
\end{corollary}

For a continuous morphism of chain complexes 
${\psi}_{\sim} : {\mathcal X}_{\sim} \rightarrow {\mathcal P}_{\sim}$
in $\fr $, a surjective map 
$ H_n(\varphi):  H_n(\X) \rightarrow H_n(\Y)$ is automatically
open, see \cite[Lemma 0.5.9]{He0}. In the category of
complete nuclear $DF$-spaces it was proved be the author in
\cite[Lemma 3.5]{Ly5}.

The following result is known for Banach and Fr\'echet spaces.

\begin{proposition}\label{coho-dual-ho}\cite[Corollary 4.9]{GLW}~
Let $(\X, d)$ be a chain complex of  Fr\'echet (Banach) spaces and
continuous linear operators
and  $ (\X^*,d^* )$ the 
strong dual cochain complex. Then the following are equivalent: 
\begin{enumerate}  
\item[(1)] $H_n(\X, d)= {\rm Ker}~ d_{n-1}/\im d_{n} $ is a Fr\'echet (Banach) space;

\item[(2)] $B_n(\X, d)= \im d_{n}$ is closed in $ X_{n}$;  
\item[(3)] $d_n$ has closed range;  
\item[(4)] the dual map $d^n=d_n^*$ has closed range;  
\item[(5)] $B^{n+1}(\X^*, d^*) = \im d^*_n$ is strongly closed in $(X_{n+1})^*$; 
\end{enumerate}
\noindent In the category of Banach spaces {\rm (1) - (5)} are equivalent to:  
\begin{enumerate}
\item[(6)] $B^{n+1}(\X^*, d^*)$ is a Banach space;  
\item[(7)] $H^{n+1}(\X^*, d^*) ={\rm Ker}~ d^*_{n+1}/\im d^*_{n}$ is a Banach space.  
\end{enumerate}

\noindent Moreover, whenever $H_n(\X, d)$ and 
$H^n(\X^*, d^*)$ are Banach spaces, up to topological isomorphism, 
$$H^n(\X^*, d^*) \iso H_n(\X, d)^*.$$
\end{proposition}

The next theorem shows that certain niceties of the theory of nuclear $DF$-spaces allow us to generalize this result to nuclear Fr\'echet spaces.

\begin{theorem}\label{H(dualX)=dualH(X)}
Let $(\X, d) $ be a chain complex of Fr\'echet spaces and continuous
linear operators and let $ (\X^*,d^* )$ be its strong dual complex.
Suppose that, for a certain $n$, either $d_n$ and  $d_{n-1}$ have 
closed ranges or $d_n^*$ and $d_{n-1}^*$ have closed ranges. 

\noindent{\rm (i)} Then, up to isomorphism of linear spaces,
$$ (H_{n}(\X, d))^* \iso H^{n}(\X^*,d^* ).$$ 

\noindent{\rm (ii)} If in addition $(\X, d) $ is a chain complex of {\em nuclear} Fr\'echet spaces, then
$H_{n}(\X, d)$ is a nuclear Fr\'echet space and, 
up to topological isomorphism,
$$ (H_{n}(\X, d))^* \iso H^{n}(\X^*,d^* )\;\;{\rm and}\;\;
H_{n}(\X, d) \iso (H^{n}(\X^*,d^*))^*.$$ 
\end{theorem}

\begin{proof} We will give a proof of (ii),  case (i) being simpler.
By \cite[Theorem 8.6.13]{Ed},  $d_n$ has closed range if and
only if $d_n^*$ has closed range. Thus  $d_{n-1}$,  $d_n$,
$d_{n-1}^*$and $d_n^*$ have closed ranges. 
We consider the following commutative diagram as in 
\cite[Lemma V.10.3]{Mac}.

\begin{equation*} 
\begin{array}{ccccccccc} ~ & ~ & 0 & ~ & 0 & ~\\ ~ & ~ & \uparrow & ~ &
\downarrow & ~\\ 0 & \rightarrow & B_n(\X) & \stackrel {\ourimath_n} 
{\hookrightarrow}
  &  Z_n(\X) & \stackrel{\sigma_n} {\rightarrow} & H_{n}(\X) & \rightarrow &
0\\ ~ & ~ & \vcenter{\llap{$\scriptstyle \tilde{d_n}~~$}}{\uparrow} & ~ &
\vcenter{\llap{$\ourjmath_n~~$}}{\downarrow} &  ~\\ ~ & ~ & X_{n+1} &
\stackrel {d_n} {\rightarrow} & X_{n} &~\\ ~ & ~ & ~ & ~ &
\vcenter{\llap{$\scriptstyle{d_{n-1}}~~$}}{\downarrow}  & ~\\ ~ & ~ & ~ & ~ & X_{n-1} & ~ \\ 
\end{array} 
\end{equation*}
where $\ourimath_n$ and $\ourjmath_n~~$ are the natural inclusions and $\sigma_n$ is the quotient map.
The notation $\tilde{d}$ is an instance of one we shall use
repeatedly, and thus we adopt the following definition.  Given a continuous linear map $\theta : E \to F$, the map $\tilde{\theta}$ is the surjective map $\tilde{\theta} : E \to \im \theta$ defined by $\tilde \theta (t) = \theta(t)$.
Here again all the maps have closed ranges.

We form the dual diagram and add the kernel of $d_n^*~$, 
$Z^{n}(\X^*) = {\rm Ker}~ d_n^*$, and the image of $d^*_{n-1}~,$
$B^n(\X^*) = \im d^*_{n-1}$ which is closed by assumption. 

\begin{equation*}
\begin{array}{ccccccccc}
 ~ & ~ & 0 & ~ & 0 & ~\\ 
 ~ & ~ & \downarrow & ~ & \uparrow & ~ \\
0 & \leftarrow &  (B_n(\X))^* & \stackrel {\ourimath_n^*}
{\longleftarrow}
  & (Z_n(\X))^* & \stackrel {\sigma_n^*} {\longleftarrow} & 
(H_{n}(\X))^* & \leftarrow & 0\\
  ~ & ~ & \vcenter{\llap{$\scriptstyle{\tilde{d}_n^*}$}}\downarrow & ~ &
\vcenter{\llap{$\scriptstyle{\ourjmath_n^*}$}}\uparrow &  ~\\ 
~ & ~ &  X_{n+1}^* & \stackrel{d_n^*} {\longleftarrow} & 
 X_{n}^* & \stackrel{\ourimath_Z} {\longleftarrow} 
& Z^n(\X^*) & \leftarrow & 0\\
 ~ & ~ & ~ & ~ & \vcenter{\llap{$\scriptstyle{d^*_{n-1}}$}}\uparrow
&~& \uparrow\vcenter{\rlap{$\scriptstyle{\ourimath_B}~~$}} & ~\\  ~ & ~ & ~
& ~ &  X_{n-1}^* & 
 \stackrel {\widetilde{d^*_{n-1}}} {\longrightarrow} & B^n(\X^*) 
& \rightarrow & 0\\ 
\end{array}
\end{equation*}
where $\widetilde{d^*_{n-1}} :  X_{n-1}^* \to \im d^*_{n-1}: \gamma
\to
[d^*_{n-1}](\gamma) =\gamma \circ d_{n-1}$. This diagram commutes and has 
exact rows and columns. 
By \cite[Lemma 2.3]{Ly4}, the exactness of
a complex  in $\fr$ is  equivalent to the exactness of its dual complex. 
Thus the exactness of the first line follows from \cite[Lemma 2.3]{Ly4}; 
of the second line from the definition of 
$Z^n(\X^*)$; of the first column  from \cite[Corollary 8.6.11]{Ed} 
since the surjectivity of $\tilde{d}_n$ implies the injectivity of $\tilde{d}_n^*$, and of the second column from \cite[Lemma 2.3]{Ly4}. Commutativity only needs to be checked for the square involving the two added terms,
namely $Z^n(\X^*)$ and $B^n(\X^*)$, and this is obvious. 
By Lemma \ref{closed_suspace_DFspaces} and Lemma \ref{KerT*barrelled},
 $Z^n(\X^*)= {\rm Ker}~  d_n^*$  and  $B^n(\X^*)= \im d^*_{n-1}$ are the strong duals of nuclear Fr\'echet spaces.
Therefore this diagram is one of strong duals of nuclear Fr\'echet spaces and continuous linear operators with closed ranges.

By Lemma \ref{KerT*barrelled}, $\im \sigma_n^*$ is a strong dual of a nuclear Fr\'echet space. Therefore, by Corollary \ref{open-map-strong-dual-FM}, the continuous linear surjective operator
$$\widetilde{\sigma_n^*}: ( H_{n}(\X))^* \to \im \sigma_n^*: \gamma
\mapsto \sigma_n^*(\gamma)$$ is open.

Let us define a map 
$$
\varphi: Z^n(\X^*) \to (H_{n}(\X))^*
$$
by the formula 
$\varphi = 
\widetilde{\sigma_n^*}^{-1} \circ \ourjmath_n^*\circ \ourimath_Z,$ where
$\widetilde{\sigma_n^*}^{-1}$ is the inverse of the topological isomorphism
$\widetilde{\sigma_n^*}$. It is now a standard diagram-chasing argument
to show that $\varphi$ is well defined and surjective.  Let us give this
argument. An element $z\in Z^n(\X^*)$ is sent by $d_n^*\circ
\ourimath_Z$ to $0$ in $X_{n+1}^*$
and therefore, since $\tilde{d_n}^*$ is injective, 
 $[\ourimath_n^* \circ \ourjmath_n^* \circ \ourimath_Z ](z) = 0$. Hence
the element $ [\ourjmath_n^*\circ \ourimath_Z](z)$ of $(Z_n(\X))^*$
belongs to ${\rm Ker}~ \ourimath_n^* = \im \sigma_n^*$,  by exactness of the
first line of the
diagram. Thus $\varphi$ is a well defined  continuous linear operator. To
show that this map is surjective, starting with $v\in (H_n(\X))^*$, we get
$u = \sigma_n^*(v) \in (Z_n(\X))^*$, 
and, since $\ourjmath_n^*$ is surjective, there is  
$t \in \X_n^*$ such that
$\ourjmath_n^*(t) = u$. It is easy to see that $t \in {\rm Ker}~ d_n^*$ 
and therefore
it lifts uniquely to $z\in Z^n(\X^*)$ and $\varphi (z) = v.$

One can see that $\ourimath_B(B^n(\X^*)) \subseteq {\rm Ker}~ \varphi, $ since
$\widetilde{d^*_{n-1}}$ is surjective and, for any 
$y \in X_{n-1}^*$, $ [\ourjmath_n^*\circ d^*_{n-1}] (y) =0.$ 
Suppose $z \in {\rm Ker}~ \varphi$, hence
$[\ourjmath_n^*\circ \ourimath_Z](z) = 0$. It implies that
$\ourimath_Z(z) \in {\rm Ker}~ \ourjmath_n^* = \im d^*_{n-1}$, so that 
there is $y  \in
X_{n-1}^*$ such that $d^*_{n-1}(y) = \ourimath_Z(z)$. Since 
$\ourimath_Z$ is
injective, $z = \ourimath_B(\widetilde{d^*_{n-1}}(y)).$ Thus ${\rm Ker}~ \varphi
=\ourimath_B[B^n(\X^*)].$ 

By Corollary \ref{open-map-strong-dual-FM}, the continuous surjective 
linear operator between  strong dual of nuclear Fr\'echet spaces
$$
\varphi: Z^n(\X^*) \to (H_{n}(\X))^*
$$
is open and, up to topological isomorphism,
$$ (H_{n}(\X, d))^* \iso H^{n}(\X^*,d^* ).$$
By \cite[Theorem 4.4.13]{Pi}, up to topological isomorphism,
$$H_{n}(\X, d) \iso (H^{n}(\X^*,d^*))^*.$$ 
\end{proof} 

\begin{corollary}\label{H(dualX) =dualH(X)-DF}
Let $(\Y, d) $ be a cochain complex of complete nuclear $DF$-spaces 
and continuous operators.
Suppose that, for a certain $n$, $d_n$ and  $d_{n-1}$ have 
closed ranges. Then
$H^{n}(\Y,d )$ is a complete nuclear $DF$-space,
$H_{n}(\Y^*, d^*)$ is a nuclear Fr\'echet space and, 
up to topological isomorphism,
$$ (H^{n}(\Y, d))^* \iso H_{n}(\Y^*,d^* )\;\;{\rm and}\;\;
H^{n}(\Y, d) \iso (H_{n}(\Y^*,d^*))^*.$$ 
\end{corollary}

\begin{proof} By \cite[Theorem 4.4.13]{Pi}, the complex $(\Y, d)$
is the strong dual of the chain complex $(\Y^*,d^* )$ of 
nuclear Fr\'{e}chet spaces and continuous linear operators.
By \cite[Theorem 4.4.12]{Pi}, complete
nuclear $DF$-spaces are reflexive, and therefore the statement
follows from Theorem \ref{H(dualX)=dualH(X)} and Proposition
\ref{coho-dual-ho}.
\end{proof}

\section{The K\"unneth formula for \frechet\ and
complete nuclear $DF$-complexes}
 
In this section we prove the existence of a topological isomorphism in the K\"unneth formula for the cohomology groups of complete nuclear $DF$-complexes (Theorem \ref{NuclearCase}). 
To start with we state the result by 
F. Gourdeau, M.C. White and the author on the K\"unneth formula for \frechet\ and Banach chain complexes. Note that similar results are true for cochain complexes. One can see that to obtain the K\"unneth formula in the category of  \frechet\   spaces  and continuous operators, we need the following notions of strict flatness (Def. \ref{strict_flat}) and of the topological purity of short exact sequences of \frechet\ spaces (Def. \ref{topol_pure}). These conditions allow us to deal with the known problems in the category of \frechet\   spaces that the projective tensor product of injective continuous linear operators is not necessarily injective and the range of an operator is not always closed.

\begin{theorem}\label{KunnethGLW} \cite[Theorem 5.2 and Corollary 4.9]{GLW}~ Let $\X$ and $\Y$ be  bounded chain
complexes in $\fr $ (in ${\mathcal Ban}$) such that all
boundary maps have closed ranges. 
Suppose that the following exact sequences of 
\frechet\ (Banach) spaces are topologically pure for all $n$:
\begin{equation}\label{Z_n(K)} 0 \to Z_n(\X) 
\stackrel{\ourimath_n}{\rightarrow} X_n
\stackrel{\tilde{d}_{n-1}}{\longrightarrow} B_{n-1}(\X)\to 0 
\end{equation}
and 
\begin{equation}\label{B_n(K)} 0 \to
B_n(\X) \stackrel{\ourjmath_n}{\rightarrow} Z_n(\X) 
\stackrel{\sigma_n} {\rightarrow}
H_{n}(\X)\to 0. 
\end{equation} 
where $\ourimath_n$ and $\ourjmath_n~~$ are the natural inclusions and $\sigma_n$ is the quotient map.
Suppose also that  $Z_n(\X)$ and $B_n(\X)$ 
are strictly flat in $\fr$ (in ${\mathcal Ban}$) for all $n$.
Then, up to  topological isomorphism, 
$$ H_n(\X \ptp\Y) \iso
\displaystyle{\bigoplus_{m+q=n}} [H_m(\X) \ptp H_q(\Y)],$$ 
and, in addition, for complexes of Banach spaces, there is also a
 topological  isomorphism 
\[ H^n((\X \ptp \Y)^*) \iso
\displaystyle{\bigoplus_{m+q=n}} [H_m(\X) \ptp H_q(\Y)]^*. \]
\end{theorem}

\begin{corollary}\label{Kunneth-Ban-Fr} Let $\X$ and $\Y$ be
bounded chain complexes  of Banach spaces and of \frechet\ spaces respectively
such that all boundary maps have closed ranges, $H_n(\X)$ and 
$B_n(\X)$ are strictly flat in ${\mathcal Ban}$ for all $n$.
Then, up to  topological  isomorphism, 
 $$ H_n(\X \ptp \Y)  \iso \displaystyle{\bigoplus_{m+q=n}} 
[H_m(\X) \ptp H_q(\Y)]$$
and, up to isomorphism of linear spaces,
$$ H^n((\X \ptp \Y)^*) \iso \displaystyle{\bigoplus_{m+q=n}} [H_m(\X) \ptp H_q(\Y)]^*.$$
If, in addition, $\Y$ is a complex  of Banach spaces, then both the above 
isomorphisms are topological.
\end{corollary}
\begin{proof} In the category of Banach spaces, by \cite[Proposition VII.1.17]{He0}, 
$B_n(\X)$ and $H_m(\X)$ strictly flat implies that $Z_n(\X)$ is strictly flat as well.
By \cite[Proposition 4.4]{Pir3}, $B_n(\X)$, $H_m(\X)$ and $Z_n(\X)$ are  also strictly flat 
in $\fr$. By \cite[Lemma 4.3]{GLW}, strict flatness of $B_n(\X)$ and $H_m(\X)$ in 
${\mathcal Ban}$ implies that the short exact sequences (\ref{Z_n(K)}) and 
(\ref{B_n(K)}) of Banach spaces are weakly admissible. By  \cite[Lemma 2.4]{Ly4},
the short exact sequences (\ref{Z_n(K)}) and  (\ref{B_n(K)}) are topologically pure in 
$\fr$. The statement follows from Theorem \ref{KunnethGLW} and  Theorem \ref{H(dualX)=dualH(X)}.

By Proposition \ref{coho-dual-ho}, in the case that both $\X$ and $\Y$ are from ${\mathcal Ban}$, we have a topological isomorphism
 $ H^n((\X \ptp \Y)^*) \iso (H_n(\X \ptp \Y))^*.$
\end{proof}

The topological isomorphism (\ref{KunnethHomology}) for homology groups under the assumptions of Part (i)  of the following theorem is already known, see M. Karoubi \cite{Ka}. To get the isomorphism for cohomology groups of dual complexes he required $H^n(\X^*)$ to be finite-dimensional.

\begin{theorem}\label{NuclearCase} Let $\X$ and $\Y$ be  bounded chain complexes in $\fr $ such that all
boundary maps have closed ranges.  

{\rm (i)} Suppose that one of complexes, say $\X$, is a complex of nuclear \frechet\ spaces.
 Then, up to  topological  isomorphism, 
\begin{equation}\label{KunnethHomology} 
H_n(\X \ptp \Y) \iso
\displaystyle{\bigoplus_{m+q=n}} [H_m(\X) \ptp H_q(\Y)]
\end{equation}
and, up to isomorphism of linear spaces,
$$ H^n((\X \ptp \Y)^*) \iso \displaystyle{\bigoplus_{m+q=n}} [H_m(\X) \ptp H_q(\Y)]^* 
\iso \displaystyle{\bigoplus_{m+q=n}} [H^m(\X^*) \ptp [H_q(\Y)]^*] .$$
{\rm (ii)} Suppose that $\X$ and $\Y$ are complexes of nuclear \frechet\ spaces.
 Then, up to topological  isomorphism,  
$$ H^n((\X \ptp \Y)^*) \iso H^n(\X^* \ptp \Y^*) \iso
\displaystyle{\bigoplus_{m+q=n}} [H^m(\X^*) \ptp H^q(\Y^*)].$$
\end{theorem}
\begin{proof} {\rm (i)} Suppose that $\X$ is a complex of nuclear \frechet\ spaces.
Since all boundary maps have closed ranges, $Z_n(\X)$ and $B_n(\X)$ are nuclear 
\frechet\ spaces. By Theorem~A.1.6 and 
Theorem~A.1.5 of \cite{EsPu},  $Z_n(\X)$ and $B_n(\X)$
are strictly flat all $n$ in  $\fr $  and the short exact sequences (\ref{Z_n(K)}) and 
 (\ref{B_n(K)}) are topologically pure in $\fr$. The first part of the statement follows from 
Theorem \ref{KunnethGLW}. By Theorem \ref{H(dualX)=dualH(X)}, up to isomorphism of linear spaces,
$ H^n((\X \ptp \Y)^*) =( H_n(\X \ptp \Y))^*$. Thus, up to isomorphism of linear spaces,
$$ H^n((\X \ptp \Y)^*) \iso \displaystyle{\bigoplus_{m+q=n}} [H_m(\X) \ptp H_q(\Y)]^* $$
By assumption, $H_m(\X)$ is a nuclear \frechet\ space for all $m$. By \cite[Theorem 21.5.9]{Ja} and by Theorem \ref{H(dualX)=dualH(X)},
up to topological isomorphism, 
$$[H_m(\X) \ptp H_q(\Y)]^*  \iso [H_m(\X)]^* \ptp [H_q(\Y)]^* \iso
H^m(\X^*) \ptp [H_q(\Y)]^*$$
for all $m, q$. 

{\rm (ii)} Since $\X$ and $\Y$ are complexes of 
nuclear \frechet\ spaces, by \cite[Theorem 21.5.9]{Ja},
up to topological isomorphism,
$(\X \ptp \Y)^* \iso \X^* \ptp \Y^*$, and so
$$ H^n((\X \ptp \Y)^*) \iso H^n(\X^* \ptp \Y^*).$$
By \cite[Proposition III.50.1]{Tr}, the projective tensor product of nuclear \frechet\ spaces is a 
nuclear \frechet\ space. Hence $\X \ptp \Y$ is a complex of 
nuclear \frechet\ spaces. By (i), for all $n$,
$$ H_n(\X \ptp\Y) \iso
\displaystyle{\bigoplus_{m+q=n}} [H_m(\X) \ptp H_q(\Y)],$$ 
is a nuclear \frechet\ space.
By Proposition \ref{coho-dual-ho} and Theorem \ref{H(dualX)=dualH(X)},
\[ H^n((\X \ptp \Y)^*) \iso (H_n(\X \ptp \Y))^* \iso 
\left(\displaystyle{\bigoplus_{m+q=n}} [H_m(\X) \ptp H_q(\Y)] \right)^*. \]
By \cite[Theorem 21.5.9]{Ja} and Theorem \ref{H(dualX)=dualH(X)}, since $H_m(\X)$ and $H_q(\Y)$ are nuclear \frechet\ spaces,
\[ 
\displaystyle{\bigoplus_{m+q=n}} [H_m(\X) \ptp H_q(\Y)]^*
\iso \displaystyle{\bigoplus_{m+q=n}} [H_m(\X)]^* \ptp [H_q(\Y)]^*
\]
\[ \iso \displaystyle{\bigoplus_{m+q=n}} [H^m(\X^*) \ptp H^q(\Y^*)].\]
\end{proof}

\begin{theorem}\label{NuclearCase-DF}~
{\rm (i)} Let  $\X$ and $\Y$ be bounded chain complexes of complete 
{\em nuclear} $DF$-spaces  such that all
boundary maps have closed ranges.  
 Then, up to  topological  isomorphism, 
 $$ H_n(\X \ptp \Y) \iso
\displaystyle{\bigoplus_{m+q=n}} [H_m(\X) \ptp H_q(\Y)].$$ 

{\rm (ii)}  Let $\X$ be a bounded chain complex of complete {\em nuclear} $DF$-spaces such that all
boundary maps have closed ranges, and let $\Y$ be a bounded chain complex of complete  $DF$-spaces such that all boundary maps of
its strong dual complex $\Y^*$ have closed ranges.   Then, up to  topological  isomorphism,  
$$ H^n((\X \ptp \Y)^*) \iso H^n(\X^* \ptp \Y^*) \iso
\displaystyle{\bigoplus_{m+q=n}} H^m(\X^*) \ptp H^q(\Y^*).$$
\end{theorem}
\begin{proof} {\rm (i)} By \cite[Theorem 4.4.13]{Pi}, the
chain complexes $\X$ and $\Y$ are the strong duals of cochain complexes $\X^*$ and $\Y^*$ of nuclear \frechet\ spaces and continuous linear operators. By Proposition \ref{coho-dual-ho}, all boundary maps of complexes $\X^*$ and $\Y^*$ have closed ranges.  By Theorem \ref{NuclearCase} (ii),
for the complexes $\X^*$ and $\Y^*$ of nuclear \frechet\ spaces,
 up to  topological  isomorphism,
$$ H_n(\X \ptp \Y) \iso H_n((\X^*)^* \ptp (\Y^*)^*) $$
$$ \iso \displaystyle{\bigoplus_{m+q=n}} [H_m((\X^*)^*) \ptp H_q((\Y^*)^*)] \iso \displaystyle{\bigoplus_{m+q=n}} [H_m(\X) \ptp H_q(\Y)].$$ 
{\rm (ii)} Since $\X$ is the  complex of complete nuclear 
$DF$-spaces, then,
by \cite[Theorem 21.5.9]{Ja}, $(\X \ptp \Y)^* \iso \X^* \ptp \Y^*$ and $$ H^n((\X \ptp \Y)^*) \iso H^n(\X^* \ptp \Y^*).$$
By Proposition \ref{coho-dual-ho}, all
boundary maps of complexes $\X^*$  have closed ranges.  
By Theorem \ref{NuclearCase} (i), for the cochain complex of nuclear \frechet\ spaces $\X^*$ (\cite[Theorem 4.4.13]{Pi}) and for the cochain complex of \frechet\ spaces $\Y^*$,
$$ H^n(\X^* \ptp \Y^*) \iso 
\displaystyle{\bigoplus_{m+q=n}} H^m(\X^*) \ptp H^q(\Y^*).$$
\end{proof}

\section{The K\"unneth formula for Hochschild cohomology of $\hat{\otimes}$-algebras which are nuclear $DF$- or \frechet\ spaces}

Let $\A$ be a $\hat{\otimes}$-algebra and let $X$ be an 
$\A$-$\hat{\otimes}$-bimodule.  We assume here that the category 
of underlying spaces $\Phi$
has the properties from Definition \ref{category}.
  Let us recall the definition of the standard homological
chain complex ${\mathcal C}_{\sim} (\A, X)$. For $n\ge 0$, let $C_n(\A, X)$
denote the projective tensor product 
 $X \ptp \A^{{\ptp}^n}$. 
The elements of $C_n(\A, X)$ are called
{\em $n$-chains}. Let the differential $d_n : C_{n+1} \to C_n$ be given by
$$ 
d_n(x \otimes  a_1 \otimes \ldots \otimes  a_ {n+1})\!= 
\! x \cdot a_1 \otimes \ldots \otimes  a_ {n+1}$$
$$ + \sum_{k=1}^{n} (-1)^k (x \otimes a_1 \otimes \ldots \otimes a_k
a_{k+1} \otimes \ldots \otimes a_ {n+1})$$
$$\;\;\;\;\;\;\; 
+(-1)^{n+1}(a_ {n+1} \cdot x  \otimes a_1 \otimes \ldots \otimes a_{n})
$$
with $d_{-1}$ the null map. The space of boundaries 
$B_n({\mathcal C}_{\sim} (\A, X))= \im d_n$ is denoted by $B_n(\A, X)$ and the space of cycles 
$Z_n({\mathcal C}_{\sim} (\A, X)) = {\rm Ker}~ d_{n-1}$ is
 denoted by $Z_n(\A, X)$. The homology groups of  this complex 
$H_n({\mathcal C}_{\sim} (\A, X)) = Z_n(\A, X)/B_n(\A, X)$ are called  
the {\it continuous Hochschild homology groups of  $\A$ with coefficients in $X$} and are 
denoted by $\H_n(\A, X)$ \cite[Definition II.5.28]{He0}.

We also consider the cohomology groups 
$H^n(({\mathcal C}_{\sim} (\A, X))^*)$ of  the dual 
complex $({\mathcal C}_{\sim} (\A, X))^*$ with the strong dual 
topology. For Banach algebras $\A$,
$({\mathcal C}_{\sim} (\A, X))^*$ is topologically isomorphic 
to the Hochschild cohomology  $\H^n(\A, X^*)$ of $\A$
with coefficients in the dual $\A$-bimodule $X^*$ 
\cite[Definition I.3.2 and Proposition II.5.27]{He0}. 

Let $\A$  be in $\Phi$  and be
a unital $\hat{\otimes}$-algebra. We put $\beta_n(\A )=
\A^{{\ptp}^{n+2}},$ $n \ge 0$ and let 
 $d_n:\beta_{n+1}(\A) \to \beta_{n}(\A)$ be given by
$$d_n( a_0 \otimes \ldots \otimes  a_ {n+2})= 
\sum_{k=0}^{n+1} (-1)^k ( a_0 \otimes \ldots \otimes a_k a_{k+1} \otimes
\ldots \otimes a_ {n+2}).$$ 
By \cite[Proposition III.2.9]{He0}, the complex over $\A$, 
$\pi:\beta(\A ) \to \A: a \otimes b \mapsto ab$,
where  $\beta(\A )$ denotes 
$$ 0 \leftarrow \beta_{0}(\A)  \stackrel{d_0}{\longleftarrow} \beta_{1}(\A) \stackrel{d_1}{\longleftarrow} \cdots 
\leftarrow \beta_n(\A )  \stackrel{d_n}{\longleftarrow} \beta_{n+1}(\A)  \leftarrow \dots $$ is a projective  resolution of 
$\A$-$\hat{\otimes}$-bimodule $\A$. $\beta(\A )$ is called 
the {\it  bar resolution} of $\A $.
The complex has a {\em contracting homotopy} 
 $s_n:\beta_n(\A )  \to \beta_{n+1}(\A) $, $(n\geq 1)$, given by 
 $$ s_n(a_0 \otimes a_1 \otimes \cdots \otimes a_{n+1}) = 1\otimes a_0 \otimes a_1 \otimes \cdots \otimes a_{n+1},  $$
 which is to say that $d_{n} s_{n} + s_{n-1} d_{n-1} = 1_{\beta_n(\A )}$.

\begin{proposition}\label{ProductResolutions}  Let $A_1$ and $A_2$ be
unital   $\hat{\otimes}$-algebras,  let $ 0 \leftarrow X 
\stackrel{\varepsilon_1}{\longleftarrow}{\mathcal X}$ 
be a projective resolution of $X \in A_1$-{\rm unmod}  and
 $ 0 \leftarrow Y \stackrel{\varepsilon_2}{\longleftarrow} {\mathcal Y}$ be a projective resolution of $Y \in A_2$-{\rm unmod}.  Then 
$ 0 \leftarrow X \hat{\otimes} Y \stackrel
{\varepsilon_1\otimes\varepsilon_2}{\longleftarrow} 
{\mathcal X} \hat{\otimes }{\mathcal Y}$
 is a projective resolution of
$X \hat{\otimes} Y \in A_1 \hat{\otimes} A_2$-{\rm unmod}. 
\end{proposition}
\begin{proof} The proof requires only minor modifications of that of 
  \cite[Proposition X.7.1]{Mac}.
\end{proof} 

  Note that the statement of Proposition \ref{ProductResolutions} 
is also true in the category of bimodules.
In the next theorem we extend the result
\cite[Theorem 6.2]{GLW} to the category of complete nuclear $DF$-spaces.

\begin{theorem}\label{ExternalProduct} Let the category for
 underlying spaces $\Phi$ be  $\fr$ or the category
of complete nuclear $DF$-spaces.
Let $\A$ and $\B$ be unital $\hat{\otimes}$-algebras with identities $e_{\A}$
and $e_{\B}$, let $X$  be an $\A$-$\hat{\otimes}$-bimodule and let $Y$ be a 
$\B$-$\hat{\otimes}$-bimodule. Then, up
to topological isomorphism, for all $n \ge 1$,
\[ \H_n(\A \ptp \B, X \ptp Y) \iso \H_n(\A \ptp \B, e_{\A} X e_{\A} \ptp e_{\B} Y e_{\B})
\]
\[  \iso H_n ({\mathcal C}_{\sim}(\A, e_{\A} X e_{\A}) \ptp 
{\mathcal C}_{\sim}(\B, e_{\B} Y e_{\B})). \]
If $X$ and $Y$ are also unital, then, up
to topological isomorphism, for all $n \ge 0$,
\[ \H_n(\A \ptp \B, X \ptp Y) \iso H_n ({\mathcal C}_{\sim}(\A, X) \ptp 
{\mathcal C}_{\sim}(\B, Y)). \]
\end{theorem}

\begin{proof} 
It is well known that, for a  $\hat{\otimes}$-algebra 
$\U$ with an identity $e$ and for a $\U$-$\hat{\otimes}$-bimodule  $Z$,
up to topological isomorphism,  for all $n \ge 1$,
$$ \H_n(\U, Z) \iso \H_n(\U, eZe), $$ where
$eZe$ is a unital $\U$-$\hat{\otimes}$-bimodule.
Therefore, up
to topological isomorphism,  for all $n \ge 1$,
\[ \H_n(\A \ptp \B, X \ptp Y) \iso \H_n(\A \ptp \B, e_{\A} X e_{\A} \ptp e_{\B} Y e_{\B}). \]

Let $\beta(\A )$ and $\beta(\B )$ be the bar
resolutions of $\A$ and $\B$.
Since the bar resolution $\beta(\A )$  is an 
 $\A$-biprojective resolution of $\A$ and  $\beta(\B )$  is a 
 $\B$-biprojective  resolution of $\B$, by
Proposition \ref{ProductResolutions} their projective tensor product
$\beta(\A )\hat{\otimes} \beta(\B )$ is an $\A\ptp\B$-biprojective resolution of $\A\hat{\otimes}\B$.

The open mapping theorem holds in the categories of \frechet\ spaces  and 
 of complete nuclear $DF$-spaces, see Corollary \ref{open-map-strong-dual-FM}   for $DF$-spaces, and,
for a continuous morphism of chain complexes 
${\psi}_{\sim} : {\mathcal X}_{\sim} \rightarrow {\mathcal P}_{\sim}$
in these  categories, a surjective map 
$ H_n(\varphi):  H_n(\X) \rightarrow H_n(\Y)$ is automatically
open, see \cite[Lemma 0.5.9]{He0} and \cite[Lemma 3.5]{Ly5}.

For a  unital $\hat{\otimes}$-algebra $\U$ and for a unital 
$\U$-$\hat{\otimes}$-bimodule  $Z$, by
\cite[Theorem III.4.25]{He0},  the Hochschild chain complex
${\mathcal C}_{\sim} (\U, Z)$ is isomorphic to 
$Z \hat{\otimes}_{\U^{e}} \beta(\U)$ and, up to topological isomorphism,
for all $n \ge 0$,
$$ \H_n(\U, Z) \iso {\rm Tor} ^{\U^e}_n(Z,\U) \iso
{\it H}_n (Z \hat{\otimes}_{\U^{e}} \beta(\U )). $$
By \cite[Section III.3.15]{He0}, the $n$th derived functor 
${\rm Tor} ^{\U^e}_n(\cdot,\U)$ does not depend on the choice of 
a $\U$-biprojective resolution of $\U$.  Therefore in these  categories, up to 
topological isomorphism,  for all $n \ge 0$,
\[
\H_n(\A \ptp \B, e_{\A} X e_{\A} \ptp e_{\B} Y e_{\B}) \iso 
{\rm Tor} ^{(\A \ptp \B)^e}_n(e_{\A} X e_{\A} \ptp e_{\B} Y e_{\B},\A \ptp \B)
\]
\[
\iso  H_n ((e_{\A} X e_{\A} \ptp e_{\B} Y e_{\B}) \hat{\otimes}_{(\A \ptp \B)^e} \beta(\A \ptp \B ))\]
\[ \iso H_n ((e_{\A} X e_{\A} \ptp e_{\B} Y e_{\B}) \hat{\otimes}_{(\A \ptp \B)^e} (\beta(\A )\hat{\otimes}\beta(\B ))).\]
By \cite[Section II.5.3]{He0}, one can prove that 
the following chain complexes are isomorphic:
\[
(e_{\A} X e_{\A} \ptp e_{\B} Y e_{\B}) \hat{\otimes}_{(\A \ptp \B)^e}
(\beta(\A )\hat{\otimes}\beta(\B )) \cong  
\]
\[(e_{\A} X e_{\A} \hat{\otimes}_{\A^e} \beta(\A ))\hat{\otimes}
(e_{\B} Y e_{\B} \hat{\otimes}_{\B^e} \beta(\B ))
\cong \]
\[ {\mathcal C}_{\sim}(\A, e_{\A} X e_{\A}) \ptp {\mathcal C}_{\sim}(\B, e_{\B} Y e_{\B})
\]
Thus, up to topological isomorphism,  for all $n \ge 0$,
\[
\H_n(\A \ptp \B, e_{\A} X e_{\A} \ptp e_{\B} Y e_{\B}) \iso 
\]
\[
 H_n ((e_{\A} X e_{\A} \ptp e_{\B} Y e_{\B}) \hat{\otimes}_{(\A \ptp \B)^e}
(\beta(\A )\hat{\otimes}\beta(\B ))) \iso 
\]
\[ H_n ({\mathcal C}_{\sim}(\A, e_{\A} X e_{\A}) \ptp {\mathcal C}_{\sim}(\B, e_{\B} Y e_{\B})).\]
\end{proof}

\begin{remark} {\rm For a  $\hat{\otimes}$-algebra 
$\U$ with an identity $e$ and for a $\U$-$\hat{\otimes}$-bimodule  $Z$,
up to topological isomorphism, for all $n \ge 1$,
$$ \H_n(\U, Z) \iso \H_n(\U, eZe), $$ where
$eZe$ is a unital $\U$-$\hat{\otimes}$-bimodule.
Thus it is easy to see that if the  boundary maps of the standard homology 
complex ${\mathcal C}_{\sim}(\U, Z)$ have closed ranges then the 
 boundary maps of the standard homology complex 
${\mathcal C}_{\sim}(\U, eZe)$ have closed ranges. The previous theorem and
this remark show that further we may concentrate on unital bimodules.}
\end{remark}

\begin{theorem}\label{Hochschild-(co)hom} Let the category for
 underlying spaces $\Phi$ be  $\fr$ or the category
of complete nuclear $DF$-spaces.
Let $\A$ and $\B$ be unital $\hat{\otimes}$-algebras, let $X$  be a unital $\A$-$\hat{\otimes}$-bimodule and let $Y$ be a unital
$\B$-$\hat{\otimes}$-bimodule. Suppose that all boundary maps of the standard homology complexes ${\mathcal C}_{\sim}(\A, X)$
and ${\mathcal C}_{\sim}(\B, Y)$ have closed ranges. Then

{\rm (i)} up to topological isomorphism in the category of complete nuclear $DF$-spaces and in the category $\fr$ under the assumption that either $\A$ and $X$ or  $\B$ and $Y$ are nuclear, for all $n \ge 0$,
\[ \H_n(\A \ptp \B, X \ptp Y) \iso
\displaystyle{\bigoplus_{m+q=n}} [\H_m(\A , X) \ptp \H_q(\B, Y)]; \]
{\rm (ii)} up to topological isomorphism in the category of complete nuclear $DF$-spaces and in the category $\fr$ under the assumption that  $\A$, $X$, $\B$ and $Y$ are nuclear, for all $n \ge 0$,
\[
H^n(({\mathcal C}_{\sim} (\A \ptp \B, X \ptp Y))^*) \iso
\left(\H_n(\A \ptp \B, X \ptp Y)\right)^* \]
\[ \iso \displaystyle{\bigoplus_{m+q=n}} 
[H^m(({\mathcal C}_{\sim} (\A, X))^*) \ptp 
H^q(({\mathcal C}_{\sim} (\B, Y))^*)]; \]
{\rm (iii)} up to isomorphism of linear spaces, in the category $\fr$ under the assumption that either $\A$ and $X$ or  $\B$ and $Y$ are nuclear, for all $n \ge 0$,
 \[
H^n(({\mathcal C}_{\sim} (\A \ptp \B, X \ptp Y))^*) \iso
\left(\H_n(\A \ptp \B, X \ptp Y)\right)^* \]
\[ \iso \displaystyle{\bigoplus_{m+q=n}} 
[H_m({\mathcal C}_{\sim} (\A, X))]^* \ptp 
[H_q({\mathcal C}_{\sim} (\B, Y))]^*. \]
\end{theorem}

\begin{proof} By Theorem \ref{ExternalProduct}, up to
topological isomorphism, for all $n \ge 0$,
\[ \H_n(\A \ptp \B, X \ptp Y)  \iso
H_n ({\mathcal C}_{\sim}(\A, X) \ptp {\mathcal C}_{\sim}(\B, Y)). \]
By \cite[Proposition III.50.1]{Tr}, the projective 
tensor product of  nuclear \frechet\ spaces is a  
nuclear \frechet\ space. By \cite[Proposition III.50.1]{Tr} and \cite[Theorem 15.6.2]{Ja}, the projective tensor product of  complete nuclear $DF$-spaces is a  complete nuclear $DF$-space. Therefore ${\mathcal C}_{\sim}(\A, X) $ and
${\mathcal C}_{\sim}(\B, Y)$ are complexes of  complete nuclear $DF$-spaces
or of [nuclear] \frechet\ spaces such that all boundary maps have closed ranges. The results follow from Theorem \ref{NuclearCase} and
Theorem \ref{NuclearCase-DF}.
\end{proof}

\begin{theorem}\label{Hochschild-(co)hom-Ban-Fr}
Let $\A$ and $\B$ be unital Banach and \frechet\  algebras respectively, let $X$  be a unital 
Banach  $\A$-bimodule and let $Y$ be a  unital  \frechet\
$\B$-bimodule. Suppose that all boundary maps of the standard homology complexes ${\mathcal C}_{\sim}(\A, X)$
and ${\mathcal C}_{\sim}(\B, Y)$ have closed ranges. Suppose that $\H_n(\A , \X)$ and $B_n(\A ,\X)$
are strictly flat in ${\mathcal Ban}$. Then,
 up to topological isomorphism, 
\[ \H_n(\A \ptp \B, X \ptp Y) \iso
\displaystyle{\bigoplus_{m+q=n}} [\H_m(\A , X) \ptp \H_q(\B, Y)], \]
and, up to isomorphism of linear spaces, 
 \[
H^n(({\mathcal C}_{\sim} (\A \ptp \B, X \ptp Y))^*) \iso
\displaystyle{\bigoplus_{m+q=n}} [\H_m(\A , X) \ptp \H_q(\B, Y)]^*. \]
\end{theorem}
\begin{proof} It follows from Theorem \ref{ExternalProduct} and Corollary \ref{Kunneth-Ban-Fr}.
\end{proof}

\begin{example} {\rm Let $\A=\ell^1(\ZZ_+)$ where  $$\ell^1(\ZZ_+)=\left\{
(a_n)_{n=0}^{\infty}: \sum_{n=0}^{\infty} |a_n| < \infty \right\}$$ be 
the unital semigroup Banach algebra of $\ZZ_+$ with  convolution
multiplication and  norm $\left\|(a_n)_{n=0}^{\infty}\right\| = 
\sum_{n=0}^{\infty} |a_n|$. In \cite[Theorem 7.4]{GLW} we showed that 
 all boundary maps of the standard homology complex ${\mathcal C}_{\sim}(\A, \A)$ have closed ranges and that $\H_n(\A , \A)$ and $B_n(\A ,\A)$
are strictly flat in ${\mathcal Ban}$. 
In \cite[Theorem 7.5]{GLW} we describe explicitly the simplicial homology groups $\H_n(\ell^1(\ZZ_+^k),\ell^1(\ZZ_+^k))$ and cohomology groups $\H^n(\ell^1(\ZZ_+^k),(\ell^1(\ZZ_+^k))^*)$ of the semigroup algebra 
$\ell^1(\ZZ_+^k)$.}
\end{example}

\begin{example} {\rm In \cite[Theorem 5.3]{Ly5} we describe explicitly the 
cyclic-type homology and cohomology groups of amenable \frechet\ algebras $\B$.
In particular we showed that all boundary maps of the standard homology complex ${\mathcal C}_{\sim}(\B, \B)$ have closed ranges. 
In \cite[Corollary 9.9]{Pir3} Pikovskii showed that an amenable unital uniform \frechet\ algebra is topologically isomorphic to the algebra $C(\Omega)$ of continuous complex-valued functions on a hemicompact $k$-space $\Omega$. Recall that a Hausdorff topological space $\Omega$ is {\em  hemicompact} if there exists a countable exhaustion $\Omega=\bigcup K_n$ with $K_n$ compact such that each compact subset of $\Omega$ is contained in some  $K_n$. A Hausdorff topological space $\Omega$ is a {\em  $k$-space} if a subset $F \subset \Omega$ is closed whenever $F \cap K$ is closed for every compact subset $K \subset \Omega$. For example, $C(\RR)$ is an amenable unital \frechet\ algebra.
}
\end{example}

The closure in a $\hat{\otimes}$-algebra ${\mathcal A}$ of the linear span of elements of the form 
$\{ab-ba:\; a, b \in {\mathcal A} \}$ is denoted by $[{\mathcal A}, {\mathcal A}].$

\begin{corollary}\label{ho-coho-l1-B} Let $\A=\ell^1(\ZZ_+)$ and let
$\I =\ell^1(\NN)$ be the closed ideal of $\ell^1(\ZZ_+)$ 
consisting of those elements with $a_0=0$. Let $C$ be an amenable 
 unital \frechet\ algebra or an amenable Banach algebra. Then
$$\H_n(\ell^1(\ZZ_+^k)\ptp \C,\ell^1(\ZZ_+^k)\ptp \C) \iso \{0 \}\;{\rm if} \;\;n>k;$$ 
$$\H^n \left( {\mathcal C}_{\sim} \left(\ell^1(\ZZ_+^k)\ptp \C,\ell^1(\ZZ_+^k) \ptp \C \right)^* \right) \iso \{0 \}\;{\rm if} \;\;n>k;$$
up to topological isomorphism, 
$$\H_n(\ell^1(\ZZ_+^k)\ptp \C,\ell^1(\ZZ_+^k)\ptp \C) \iso 
{\bigoplus\nolimits^{k \choose n}} \left(\I^{\ptp^n} \ptp
\A^{\ptp^{k-n}} \right) \ptp \left(\C/[\C, \C] \right)$$
if $ n\le k;$
\noindent and, up to isomorphism of linear spaces for \frechet\  algebras $\C$
and up to topological isomorphism for Banach algebras $\C$,
$$\H^n \left( {\mathcal C}_{\sim} \left(\ell^1(\ZZ_+^k)\ptp \C,\ell^1(\ZZ_+^k) \ptp \C \right)^* \right) \iso {\bigoplus\nolimits^{k \choose n}}  \left(\I^{\ptp^n} \ptp \A^{\ptp^{k-n}}\ptp \left(\C/[\C,\C]\right) \right)^* $$
if $n\le k.$
Moreover, for Banach algebras $\C$, up to topological isomorphism, for all $n \ge 0$,
$$\H^n(\ell^1(\ZZ_+^k)\ptp \C,(\ell^1(\ZZ_+^k)\ptp \C)^*) \iso \H^n({\mathcal C}_{\sim}(\ell^1(\ZZ_+^k)\ptp \C,\ell^1(\ZZ_+^k)\ptp \C)^*).$$
\end{corollary}

\begin{proof} By \cite[Theorem 5.3]{Ly5}, for an amenable \frechet\ algebra $\C$, \newline ${\mathcal H}_0({\mathcal C}, \C) \iso {\mathcal C}/[{\mathcal C}, {\mathcal C}]$ and ${\mathcal H}_{n}({\mathcal C}, \C) \iso \{0 \}$ for  all  $n \ge 1$. Recall that an amenable Banach algebra has a bounded approximate identity.

In \cite[Theorem 7.4]{GLW} we showed that 
 all boundary maps of the standard homology complex ${\mathcal C}_{\sim}(\A, \A)$ have closed ranges and that $\H_n(\A , \A)$ and $B_n(\A ,\A)$
are strictly flat in ${\mathcal Ban}$. By \cite[Proposition 7.3]{GLW},
up to topological isomorphism, the
simplicial homology groups $\H_n(\A,\A)$ are given by 
  $\H_0(\A,\A) \iso \A = \ell^1(\ZZ_+),$
  $\H_1(\A,\A) \iso \I = \ell^1(\NN),$
  $\H_n(\A,\A) \iso \{0\} $ for $n\geq 2$.\\

Note that $\ell^1(\ZZ_+^k)\ptp \C \iso \A \ptp \B$ where $\B=\ell^1(\ZZ_+^{k-1})\ptp \C$. We use induction on $k$ to prove the corollary for homology groups. 
For $k=1$, the result  follows from Theorem \ref{Hochschild-(co)hom-Ban-Fr}  for an amenable 
 unital \frechet\ algebra $\C$, and from \cite[Theorem 5.5]{GLW2} for an amenable Banach algebra $\C$. The
simplicial homology groups $\H_n(\A \ptp \C,\A \ptp \C)$ are given,
up to topological isomorphism, by \begin{eqnarray*}
  \H_0(\A \ptp \C,\A \ptp \C) & \iso & \A \ptp \left(\C/[\C, \C] \right), \\
  \H_1(\A \ptp \C,\A \ptp \C) & \iso & \I \ptp \left(\C/[\C, \C] \right), \\
  \H_n(\A \ptp \C,\A \ptp \C) &\iso \{0\} & \hbox{ for $n\geq 2$} .
\end{eqnarray*}
Let $k>1$ and 
suppose that the result for homology holds for $k-1$. 
As  $\ell^1(\ZZ_+^k)\ptp \C \iso \A \ptp \B$ where $\B=\ell^1(\ZZ_+^{k-1})\ptp \C$, we have
$$\H_n(\ell^1(\ZZ_+^k)\ptp \C,\ell^1(\ZZ_+^k)\ptp \C) \iso \H_n(\A \ptp \B, \A \ptp \B).$$ 
Also, it follows from the inductive hypothesis that, 
for all $n$, the $\H_n(\B,\B)$ are \frechet\ [Banach] spaces and hence the $B_n(\B,\B)$ are closed. 
We can therefore apply Theorem \ref{Hochschild-(co)hom-Ban-Fr} for an amenable 
 unital \frechet\ algebra $\C$ and \cite[Theorem 5.5]{GLW2} for an amenable Banach algebra $\C$, to get
\begin{eqnarray*}  \H_n(\A \ptp \B, \A \ptp \B) & \iso & 
 \bigoplus_{m+q=n}\left[\H_m(\A,\A)\ptp \H_q(\B,\B)\right]. \end{eqnarray*}
The terms in this direct sum vanish for $m\ge 2$, and thus we only need 
to consider 
$$\left(\H_0(\A,\A)\ptp \H_n(\B,\B)\right) \oplus 
\left(\H_1(\A,\A)\ptp \H_{n-1}(\B,\B)\right).$$ 
The rest is clear.
\end{proof}

\begin{example} {\rm 
Let $\A$ be  the convolution Banach algebra $L^1(\RR_+)$ of
complex-valued, Lebesgue measurable functions $f$ on $\RR_+$ with
finite $L^1$-norm. 
In \cite[Theorem 4.6]{GLW2} we showed that 
 all boundary maps of the standard homology complex ${\mathcal C}_{\sim}(\A, \A)$ have closed ranges and that $\H_n(\A , \A)$ and $B_n(\A ,\A)$
are strictly flat in ${\mathcal Ban}$. 
In \cite[Theorem 6.4]{GLW2} we describe explicitly the simplicial homology groups $\H_n(L^1(\RR_+^k),L^1(\RR_+^k))$ and cohomology groups $\H^n(L^1(\RR_+^k),(L^1(\RR_+^k))^*)$ of the semigroup algebra 
$L^1(\RR_+^k)$.}
\end{example}

\begin{corollary}\label{ho-coho-L1-B} Let $\A=L^1(\RR_+^k)$ and let
 $C$ be an amenable Banach algebra. Then
$$\H_n( L^1(\RR_+^k)\ptp \C, L^1(\RR_+^k)\ptp \C) \iso \{0 \}\;{\rm if} \;\;n>k;$$ 
$$\H^n \left(L^1(\RR_+^k)\ptp \C, \left(L^1(\RR_+^k)\ptp \C \right)^* \right) 
\iso \{0 \}\;{\rm if} \;\;n>k;$$
up to topological isomorphism, 
$$\H_n(L^1(\RR_+^k)\ptp \C, L^1(\RR_+^k)\ptp \C) \iso {\bigoplus\nolimits^{k \choose n}} L^1(\RR_+^k)\ptp \left(\C/[\C, \C] \right)\;{\rm if} \; n\le k;$$
\noindent and
$$\H^n(L^\infty(\RR_+^k)\ptp \C, (L^\infty(\RR_+^k)\ptp \C)^*) \iso
 {\bigoplus\nolimits^{k \choose n}} \left[L^\infty(\RR_+^k) \ptp \left(\C/[\C,\C] \right) \right]^* $$
if $n\le k.$
\end{corollary}

\begin{proof} Note that $\A$ and $\C$ have bounded approximate identities.
By \cite[Theorem 5.3]{Ly5}, for an amenable Banach algebra $\C$, 
${\mathcal H}_{n}({\mathcal C}, \C) \iso \{0 \}$ for  all  $n \ge 1$,
${\mathcal H}_0({\mathcal C}, \C) \iso {\mathcal C}/[{\mathcal C}, {\mathcal C}]$. Therefore all boundary maps of the standard homology complex ${\mathcal C}_{\sim}(\C, \C)$ have closed ranges.

In \cite[Theorem 4.6]{GLW2} we showed that 
 all boundary maps of the standard homology complex ${\mathcal C}_{\sim}(\A, \A)$ have closed ranges and that $\H_n(\A , \A)$ and $B_n(\A ,\A)$
are strictly flat in ${\mathcal Ban}$. By \cite[Theorem 4.6]{GLW2},
up to topological isomorphism, the
simplicial homology groups $\H_n(\A,\A)$ are given by 
  $\H_0(\A,\A) \iso \H_1(\A,\A) \iso \A = L^1(\RR_+)$ and
  $\H_n(\A,\A) \iso \{0\} $ for $n\geq 2$.\\

 Note that 
$L^1(\RR_+^k)\ptp \C \iso \A \ptp \B$ where $\B=  L^1(\RR_+^{k-1}) \ptp \C$. We use induction on $k$ to prove the corollary for homology groups. 
For $k=1$, the result follows  from \cite[Theorem 5.5]{GLW2}. The simplicial homology groups $\H_n(\A \ptp \C,\A \ptp \C)$ are given, up to topological isomorphism, by 
 $$ \H_0(\A \ptp \C,\A \ptp \C) \iso  \H_1(\A \ptp \C,\A \ptp \C) \iso \A \ptp \left(\C/[\C, \C] \right)$$ and
 $$ \H_n(\A \ptp \C,\A \ptp \C) \iso \{0\}$$  for $n\geq 2$.

Let $k>1$ and suppose that the result for homology holds for $k-1$. 
As $ L^1(\RR_+^k)\ptp \C \iso \A \ptp \B$ where $\B=  L^1(\RR_+^{k-1}) \ptp \C$, we have
$$\H_n(L^1(\RR_+^k)\ptp \C,L^1(\RR_+^k)\ptp \C) \iso \H_n(\A \ptp \B, \A \ptp \B).$$ 
Also, it follows from the inductive hypothesis that, 
for all $n$, the $\H_n(\B,\B)$ are Banach spaces and hence the $B_n(\B,\B)$ are closed. 
We can therefore apply \cite[Theorem 5.5]{GLW2} for an amenable Banach algebra $\C$, to get
\begin{eqnarray*}  \H_n(\A \ptp \B, \A \ptp \B) & \iso & 
 \bigoplus_{m+q=n}\left[\H_m(\A,\A)\ptp \H_q(\B,\B)\right]. \end{eqnarray*}
The terms in this direct sum vanish for $m\ge 2$, and thus we only need 
to consider 
$$\left(\H_0(\A,\A)\ptp \H_n(\B,\B)\right) \oplus 
\left(\H_1(\A,\A)\ptp \H_{n-1}(\B,\B)\right).$$

For cohomology groups, by \cite[Corollary 4.9]{GLW},
$$\H^n(L^\infty(\RR_+^k)\ptp \C, (L^\infty(\RR_+^k)\ptp \C)^*) \iso 
\H^n \left({\mathcal C}_{\sim} \left(L^1(\RR_+^k)\ptp \C, L^1(\RR_+^k) \ptp \C \right)^* \right)$$
$$ \iso {\bigoplus\nolimits^{k \choose n}} \left[L^\infty(\RR_+^k) \ptp \left(\C/[\C,\C] \right) \right]^* $$
if $n\le k.$
\end{proof}

\section{Applications to the cyclic-type cohomology of 
certain  \frechet\ and $DF$ algebras}

In this section we give explicit formulae for the continuous cyclic-type homology and cohomology of projective tensor products of certain $\hat{\otimes}$-algebras which are \frechet\ spaces or complete nuclear $DF$-spaces.

One can consult the books by Loday \cite{Loday} or Connes 
\cite{Co2} on cyclic-type homological theory.
The  {\em continuous bar} and {\em \lq naive' Hochschild
homology of} a  $\hat{\otimes}$-algebra $\A$ are 
defined respectively as 
$$
{\mathcal H}^{bar}_*(\A) = H_*({\mathcal C}(\A), b') 
\;\; {\rm and}\;\; 
{\mathcal H}^{naive}_*(\A)= H_*({\mathcal C} (\A), b),$$
where  ${\mathcal C}_n(\A) =  \A^{\hat{\otimes}(n+1)}$,  and the differentials
$b$, $b'$ are given by 
$$
b'(a_0 \otimes \dots \otimes a_n) = \sum_{i=0}^{n-1} (-1)^i (a_0 \otimes \dots \otimes a_i a_{i+1}\otimes \dots \otimes a_n) \; {\rm and} \;
$$
$$
b(a_0 \otimes \dots \otimes a_n) = 
b'(a_0 \otimes \dots \otimes a_n) + 
(-1)^n(a_n a_0 \otimes \dots \otimes a_{n-1}).
$$
Note that ${\mathcal H}^{naive}_*(\A)$ is just another way of writing ${\mathcal H}_*(\A, \A)$,
the continuous  homology of $\A$ with coefficients in $\A,$ as described in \cite{He0,BEJ1}. 

For a $\hat{\otimes}$-algebra $\A$, consider the mixed complex 
$(\bar{\Omega} \A_+, \tilde{b}, \tilde{B})$, where
$\bar{\Omega}^n \A_+ = \A^{\hat{\otimes}(n+1)} \oplus 
\A^{\hat{\otimes} n}$ and
$$
\tilde{b} = \larray b & 1-\lambda \\ 0 & -b'\rarray; \;\;\;
\tilde{B} = \larray 0 & 0 \\ N & 0 \rarray 
$$
where $\lambda(a_1 \otimes \dots \otimes a_{n}) = 
(-1)^{n-1}(a_n \otimes a_1 \otimes \dots \otimes a_{n-1})$ and
$ N = {\rm id} + \lambda + \dots + \lambda^{n-1}$  \cite[1.4.5]{Loday}.  
The {\it continuous Hochschild homology of $\A$}, the {\it continuous cyclic
homology of $\A$ }  and the {\it  continuous 
periodic cyclic homology of $\A$ } are
defined by 
$$
{\mathcal H}{\mathcal H}_*(\A) = H^b_*(\bar{\Omega} \A_+, \tilde{b}, \tilde{B}),
\; \;
{\mathcal H}{\mathcal C}_*(\A) = H^c_*(\bar{\Omega} \A_+, \tilde{b}, \tilde{B})
\; \;{\rm and} $$
$$
{\mathcal H}{\mathcal P}_*(\A) = H^p_*(\bar{\Omega} \A_+, \tilde{b}, \tilde{B})
$$
where $H^b_*$, $H^c_*$ and $H^p_*$ are 
Hochschild homology, cyclic homology and periodic cyclic homology
of the mixed complex $(\bar{\Omega} \A_+, \tilde{b}, \tilde{B})$ in
the category $\mathcal{LCS}$ of locally convex spaces and 
continuous linear operators; see, for example, \cite{Ly3}.

There is also a {\it cyclic  cohomology} theory associated 
with a complete locally convex algebra $\A,$
obtained when one replaces the chain complexes of $\A$ by their dual 
complexes of strong dual spaces. For example, the {\em continuous
bar cohomology} ${\mathcal H}^n_{bar}(\A)$ of $\A$ is the cohomology of the dual complex $({\mathcal C}(\A)^*, (b')^*)$ of 
$({\mathcal C}(\A), b')$. \\

\begin{example}~  {\rm  Some examples of $C^*$-algebras without
non-zero bounded traces are: \\
(i)  The $C^*$-algebra ${\cal K}(H)$ of compact operators on  an infinite-dimensional Hilbert space $H$; see \cite[Theorem 2]{An}. We can also show that \newline $ C(\Omega,{\cal K}(H))^{tr} = 0$, where $\Omega$ is a compact space. \\
 (ii)  Properly infinite von Neumann
algebras ${\cal U}$; see \cite[Example 4.6]{Ly2}. This class includes the   $C^*$-algebra ${\cal B}(H)$ of all bounded operators on  an infinite-dimensional Hilbert space $H$; see also \cite{Hal} for the statement ${\cal B}(H)^{tr} = 0$.}\\
\end{example}

\begin{example}~ {\rm   Let $H = 
\displaystyle{ \lim_{\stackrel{\textstyle\rightarrow}{  i}}}~
 H_i$ be a  strict inductive limit of Hilbert spaces.
 Suppose that
$H_1$ and $H_{m+1}/H_{m}, \; m=1,2, \dots$, are
 infinite-dimensional  spaces. 
Consider the Fr\'echet locally $C^*$-algebra
$ \mathcal{L}(H)$ of  continuous linear operators $T$ on $H$
that leave each $ H_i$ invariant and satisfy
 $T_j P_{ij}= P_{ij}T_j $ for all $i<j$
where $T_j= T|H_j: T_j(\eta) = T(\eta)$ for $\eta \in H_j$ and 
$P_{ij}$ is the projection from $ H_j$ onto $ H_i$.
By \cite[Example 6.6]{Ly3} that, for all $n\ge 0$,
$\H_n(\mathcal{L}(H),\mathcal{L}(H))=\{0\}$.
}
\end{example}

\begin{corollary}\label{ho-coho-A-B(H)}
Let ${\mathcal A}$ be a \frechet\ algebra belonging to
one of the following classes:

{\rm (i)} $\A=\ell^1(\ZZ_+^k) \ptp \C$, where 
 $\C$ is a $C^*$-algebra without non-zero bounded traces;

{\rm (ii)} $\A=L^1(\RR_+^k) \ptp \C$, where 
 $\C$ is a $C^*$-algebra without non-zero bounded traces;

{\rm (iii)} $\A=\ell^1(\ZZ_+^k) \ptp \C$, where 
 $\C$ is  the \frechet\ locally $C^*$-algebra ${\mathcal L}(H)$ of continuous linear operators on a strict inductive limit 
 $H = 
\displaystyle{ \lim_{\stackrel{\textstyle\rightarrow}{  i}}}~H_i$
 of Hilbert spaces
such that $H_1$ and $H_{m+1}/H_{m}, \; m=1,2, \dots$, are
 infinite-dimensional  spaces; 

{\rm (iv)} $\A=\D \ptp \C$, where $\D$ is a unital nuclear \frechet\ algebra such that all boundary maps of the standard homology complex ${\mathcal C}_{\sim}(\D, \D)$ have closed ranges
 and $\C$ is  the \frechet\ locally $C^*$-algebra ${\mathcal L}(H)$ of continuous linear operators on a strict inductive limit 
 $H = 
\displaystyle{ \lim_{\stackrel{\textstyle\rightarrow}{  i}}}~ H_i$ of Hilbert spaces
such that $H_1$ and $H_{m+1}/H_{m}, \; m=1,2, \dots$, are
 infinite-dimensional  spaces. 
 
Then, $ \H_n(\A,  \A)  \iso  \{0\}$
and $ \H^n(\A,  \A) \iso \{0\}$ for all $n \ge 0$;

$$HH_n(\A) \iso HH^n(\A) \iso \{0\}\;\;{\rm for\;\; all}\;\; n\ge 0,$$
$$HC_n(\A) \iso  HC^n(\A) \iso \{0\}\;\;{\rm for\;\; all}\;\; n\ge 0,$$
and
$$HP_m(\A) \iso HP^m(\A) \iso \{0\}\;\;{\rm for}\;\; m=0,1.$$
\end{corollary}

\begin{proof}  
By \cite[Theorem 4.1 and Corollary 3.3]{ChSi}, for a
$C^*$-algebra ${\cal U}$ without non-zero bounded traces $\H^n({\cal U}, {\cal U}^*) \iso  \{0\}$ for all   $n \ge 0$. 
By \cite[Example 6.6]{Ly3}, for the \frechet\ locally $C^*$-algebra 
${\mathcal L}(H)$, 
\newline $\H_n( {\mathcal L}(H),{\mathcal L}(H)) \iso  \{0\}$ for all   $n \ge 0$. 
 
In cases (i) and (ii) we use induction by $k$ and apply  \cite[Theorem 5.5]{GLW2} for a $C^*$-algebra $\C$ without non-zero bounded traces, to get
$ \H_n(\A,  \A) \iso \{0\}$ for all $n \ge 0$. For example,
as in Corollary \ref{ho-coho-l1-B}, note that $\ell^1(\ZZ_+^k)\ptp \C \iso \ell^1(\ZZ_+) \ptp \B$ where $\B=\ell^1(\ZZ_+^{k-1})\ptp \C$. The algebras $\ell^1(\ZZ_+)$ and $\B$ satisfy the conditions of Theorem \ref{Hochschild-(co)hom-Ban-Fr}. We use induction on $k$ to prove that 
$$ \H_n(\ell^1(\ZZ_+^k) \ptp \C,  \ell^1(\ZZ_+^k) \ptp \C) \iso \{0\}$$
for all $n \ge 0$. 

In case (iii) we use induction on $k$ and apply  Theorem \ref{Hochschild-(co)hom-Ban-Fr} for the \frechet\ locally $C^*$-algebra ${\mathcal L}(H)$, to get
$ \H_n(\A,  \A) \iso \{0\}$ for all $n \ge 0$. 

In case (iv) we apply Theorem \ref{Hochschild-(co)hom}, to get
$ \H_n(\A,  \A) \iso \{0\}$ for all $n \ge 0$. 

The triviality of the continuous cyclic and periodic cyclic homology and cohomology groups follows from \cite[Corollory 4.7]{Ly3}.
\end{proof}

 A  $\hat{\otimes}$-algebra $\A$ is said to be 
{\it biprojective} if it is projective in
the category of  $\A$-$\hat{\otimes}$-bimodules  \cite[Def. 4.5.1]{He0}. 
A  $\hat{\otimes}$-algebra $\A$ is said to be 
{\it contractible} if $\A_+$ is projective 
  $\A$-$\hat{\otimes}$-bimodules. A  $\hat{\otimes}$-algebra $\A$ is
contractible if and only if $\A$ is biprojective and has an identity
 \cite[Def. 4.5.8]{He0}.

Recall that, for a $\hat{\otimes}$-algebra $\A$ and for an 
$A$-$\hat{\otimes}$-bimodule $X$, 
 $[X, \A]$ is the closure in $X$ of the linear span of elements of the
form $ a\cdot x - x \cdot a; x \in X, a \in \A$;
$${\rm Cen}_{\A} X
 = \{ x \in X : a\cdot x = x \cdot a \;
{\rm for} \;{\rm all} \;  a \in \A \} \;\; {\rm and}$$
$${\rm Cen}_{\A} X^* = \{ f \in X^* : f( a\cdot x) = f(x \cdot a) \;
{\rm for} \;{\rm all} \;  a \in \A \;{\rm and} \; x \in X. \}$$

\begin{theorem}\label{A-X-trivial-Kunn-hom} Let the category for
 underlying spaces $\Phi$ be  $\fr$ or the category
of complete nuclear $DF$-spaces. Let $\A$ and $\B$ be unital $\hat{\otimes}$-algebras, let $Y$ be a unital $\B$-$\hat{\otimes}$-bimodule and let $X$  be a unital $\A$-$\hat{\otimes}$-bimodule such that ${\mathcal H}_0( \A, X)$ is Hausdorff and ${\mathcal H}_n( \A, X) = \{0 \}$ for  all  $n \ge 1$. Suppose that all boundary maps of the standard homology complex ${\mathcal C}_{\sim}(\B, Y)$ have closed ranges. Then

{\rm (i)} up to topological isomorphism in the category of complete nuclear $DF$-spaces and in the category $\fr$ under the assumption that either $\A$ and $X$ or  $\B$ and $Y$ are nuclear, for all $n \ge 0$,
\[ \H_n(\A \ptp \B, \X \ptp Y) \iso   X/[X, \A] \ptp \H_n(\B, Y);\]

{\rm (ii)} up to topological isomorphism in the category of complete nuclear $DF$-spaces and in the category $\fr$ under the assumption that  $\A$, $X$, $\B$ and $Y$ are nuclear, for all $n \ge 0$,
\[
H^n(({\mathcal C}_{\sim} (\A \ptp \B, X \ptp Y))^*) \iso
\left(\H_n(\A \ptp \B, X \ptp Y)\right)^* \]
\[ \iso {\rm Cen}_{\A} X^* \ptp 
H^n(({\mathcal C}_{\sim} (\B, Y))^*); \]

{\rm (iii)} up to isomorphism of linear spaces, in the category $\fr$ under the assumption that either $\A$ and $X$ or  $\B$ and $Y$ are nuclear, for all $n \ge 0$,
\[
H^n(({\mathcal C}_{\sim} (\A \ptp \B, X \ptp Y))^*) \iso 
\left(\H_n(\A \ptp \B, X \ptp Y)\right)^* \]
\[ \iso (X/[X, \A] )^* \ptp (\H_n(\B, Y))^*. \]
\end{theorem}

\begin{proof} By assumption ${\mathcal H}_0( \A, X)$ is Hausdorff, 
and so  ${\mathcal H}_0( \A, X) \iso X/[X, \A]$.
The result follows from Theorem \ref{Hochschild-(co)hom}.
\end{proof}

\begin{lemma} Let the category for underlying spaces $\Phi$ be  $\fr$ or the category of complete nuclear $DF$-spaces.
Let  ${\mathcal A}$ be  a contractible $\hat{\otimes}$-algebra.  Then, for each $\A$-$\hat{\otimes}$-bimodule $X$,   ${\mathcal H}_0( \A, X)$ is Hausdorff and ${\mathcal H}_n( \A, X) \iso \{0 \}$ for  all  $n \ge 1$. 
\end{lemma}

\begin{proof} In the categories of \frechet\ spaces and complete nuclear $DF$-spaces, the open mapping theorem holds. Therefore, by \cite[Theorem III.4.25]{He0}, for all $n \ge 0$ and all $\A$-$\hat{\otimes}$-bimodule $X$, up to topological isomorphism,
$$ \H_n(\A,X) \iso {\rm Tor}_n^{\A^e}(X,\A_+).$$
Since  $\A$ is contractible, by \cite[Proposition III.3.5 and Proposition VII.1.2]{He0}, $\H_n(\A,X)  \iso \{0\}$ for all $n\ge 1$ and
$ \H_0(\A,X) \iso {\rm Tor}_0^{\A^e}(X,\A_+)$
is Hausdorff.
\end{proof}

\subsection{Contractible $\hat{\otimes}$-algebras}

\begin{example}~ 
{\rm A countable direct product $\prod_{i \in J}  M_{n_i}(\CC)$ of full matrix algebras 
is contractible \frechet\ algebra \cite{Tay1}.}\\
\end{example}

\begin{example}~ 
{\rm Let  $G$ be a compact Lie group and let 
 $\mathcal{E}^*(G)$ be the strong dual to the nuclear \frechet\ algebra of smooth functions $\mathcal{E}(G)$ on $G$ with the convolution product, so that  $\mathcal{E}^*(G)$ is a complete nuclear 
$DF$-space. This is a $\hat{\otimes}$-algebra with respect to convolution 
multiplication: 
for $f,g \in \mathcal{E}^*(G)$ and $x \in  \mathcal{E}(G)$, 
$ <f*g, x>= <f, y>$, where $y \in \mathcal{E}(G)$ is defined by 
$ y(s)= <g, x_s>, \; s \in G$ and $x_s(t)= x (s^{-1}t),\; t \in G$.
J.L. Taylor proved that the algebra of distributions
$\mathcal{E}^*(G)$ on a compact Lie group $G$ is contractible 
 \cite{Tay1}.\\
}
\end{example}

\begin{example}~{\rm  Fix a real number $1 \le R \le \infty$ and a
nondecreasing sequence $\alpha = (\alpha_i)$ of positive numbers with
$\lim_{i \to  \infty} \alpha_i = \infty$. The power series space
$$ \Lambda_R(\alpha)= 
\{ x = (x_n) \in \CC^\NN : \|x \|_r = \sum_n |x_n| r^{\alpha_n}
< \infty \;{\rm for} \;{\rm all} \; 0 <r <R \}$$
is a \frechet\  K\"{o}the algebra with pointwise multiplication.
The topology  of $ \Lambda_R(\alpha)$ is  determined by
a countable  family of seminorms  $\{\|x \|_{r_k} :k \in \NN\}$ 
where $\{r_k\}$ is an arbitrary increasing sequence converging to $R$.

By \cite[Corollary 3.3]{Pir2}, $ \Lambda_R(\alpha)$ is {\it biprojective }
 if and only if $R=1$ or $R=\infty$.

By the Grothendieck-Pietsch criterion, $ \Lambda_R(\alpha)$ is 
{\it nuclear} if and only if  for $ \overline{\lim_n }\frac{\log n}{\alpha_n} = 0$ for
$R < \infty$ and $ \overline{\lim_n }\frac{\log n}{\alpha_n}  < \infty$
for $R = \infty$, see  \cite[Example 3.4]{Pir1}.

By \cite[Proposition 3.15]{Pir2}, for the \frechet\  K\"{o}the 
algebra $ \Lambda_1(\alpha)$, the following conditions are
equivalent:
(i)  $ \Lambda_1(\alpha)$ is  contractible,
(ii) $ \Lambda_1(\alpha)$  is nuclear,
(iii) $ \Lambda_1(\alpha)$  is unital.

By \cite[Corollary 3.18]{Pir2}, if 
 $ \Lambda_{\infty}(\alpha)$ is nuclear,  then
 the strong dual 
$\Lambda_{\infty}(\alpha)^*$ is a nuclear, contractible K\"{o}the 
$\hat{\otimes}$-algebra which is a  $DF$-space.
}
\end{example}

The algebra $ \Lambda_R((n))$ is topologically isomorphic to the 
algebra of functions holomorphic on the open disc of radius $R$,
 endowed with {\it Hadamard product}, that is, with 
``co-ordinatewise" product of the Taylor expansions of holomorphic 
functions.\\

\begin{example}~ {\rm
The algebra $\mathcal{H}(\CC) \iso  \Lambda_{\infty}((n))$ 
of entire functions, endowed with the Hadamard product, is a
biprojective nuclear \frechet\ algebra \cite{Pir2}. The strong dual 
$\mathcal{H}(\CC)^*$ is a nuclear contractible K\"{o}the 
$\hat{\otimes}$-algebra which is a  $DF$-space.}\\
\end{example}

\begin{example}~{\rm
The algebra $\mathcal{H}({\DD}_1) \iso \Lambda_1((n))$ of 
functions holomorphic on the open unit disc,
 endowed with the Hadamard product, is a
biprojective nuclear \frechet\ algebra. Moreover it is contractible,
since the function $z \mapsto (1- z)^{-1}$ is an identity for
 $\mathcal{H}({\DD}_1)$ \cite{Pir2}.}\\
\end{example}

\begin{example}\label{s-algebra}~{\rm
 The nuclear \frechet\ algebra of rapidly decreasing
sequences 
 $$s = \{ x = (x_n) \in \CC^\NN : \|x \|_k = \sum_n |x_n| n^k
< \infty \;{\rm for} \;{\rm all} \; k \in \NN \}$$
is a biprojective  K\"{o}the algebra \cite{Pir1}. The algebra $s$
is topologically isomorphic to $ \Lambda_{\infty}(\alpha)$
with $\alpha_n =\log n$ \cite{Pir2}. The nuclear K\"{o}the 
$\hat{\otimes}$-algebra $s^*$ 
of sequences of polynomial growth is contractible \cite{Tay1}.}\\
\end{example} 

The space of continuous traces on a topological algebra $\A$ is denoted by $\A^{tr}$, that is,
$$\A^{tr} = {\rm Cen}_{\A} \A^*= \{f \in \A^*: f(ab) = f(ba)\;{\rm for\; all}\;a, b \in \A \}.$$
The closure in $\A$ of the linear span of elements of the form 
$\{ab-ba:\; a, b \in \A \}$ is denoted by $[\A, \A].$
Recall that $b_0: \A \hat{\otimes} \A \to \A$ is uniquely 
determined by $a \otimes b \mapsto ab -ba$.\\

\begin{corollary}\label{AandB-simpl-trivial-Kunn-hom}~ Let the category for
 underlying spaces $\Phi$ be  $\fr$ or the category
of complete nuclear $DF$-spaces. Let $\A$ and $\B$ be unital $\hat{\otimes}$-algebras such that ${\mathcal H}_0( \A, \A)$ and ${\mathcal H}_0( \B, \B)$ are
Hausdorff, and $${\mathcal H}_n( \A, \A) \iso {\mathcal H}_n( \B, \B) \iso
\{0 \}$$
for  all  $n \ge 1$. Then

{\rm (i)}~ up to topological isomorphism in the category of complete nuclear $DF$-spaces and in the category $\fr$ under the assumption that either $\A$  or  $\B$ is nuclear, 
\begin{equation}
\begin{array}{ccccccccccc} \label{nice-homology-AB}
 {\mathcal H}{\mathcal H}_{0}(\A\ptp \B) \iso {\mathcal H}^{naive}_{0}(\A \ptp \B) \iso \A/[\A, \A]\ptp \B/[\B,\B]\;{\rm and} \;\\
{\mathcal H}{\mathcal H}_{n}(\A\ptp \B) \iso {\mathcal H}^{naive}_n
(\A \ptp \B) \iso \{0 \} \;\;{\rm for} \; {\rm  all}\;n \ge 1;\\
&\\
{\mathcal H}{\mathcal C}_{2\ell}(\A\ptp \B) \iso
\A/[\A, \A] \ptp \B/[\B,\B] \;{\rm and} \;\\
{\mathcal H}{\mathcal C}_{2\ell+1}(\A\ptp \B) \iso \{0\} \;{\rm for} \; {\rm  all}\; \ell \ge 0; \\
\end{array}
\end{equation}

{\rm (ii)}~ up to topological isomorphism in the category of complete nuclear $DF$-spaces and in the category $\fr$ under the assumption that  $\A$ and $\B$  are nuclear,
\begin{equation}
\begin{array}{ccccccccccc} \label{nice-cohomology-AB}
{\mathcal H}{\mathcal H}^{0}(\A\ptp \B) \iso
{\mathcal H}_{naive}^{0}(\A\ptp \B) \iso (\A\ptp \B)^{tr}\;{\rm and} \;\\
{\mathcal H}{\mathcal H}^{n}(\A\ptp \B) \iso {\mathcal H}_{naive}^n(\A\ptp \B) \iso \{0 \} \;\;{\rm for} \; {\rm  all}\;
 n \ge 1;\\ 
&\\
{\mathcal H}{\mathcal C}^{2\ell}(\A\ptp \B) \iso (\A\ptp \B)^{tr} \;{\rm and} \;\\
{\mathcal H}{\mathcal C}^{2\ell+1}(\A\ptp \B) \iso \{0\}  
\;\;{\rm for} \; {\rm  all}\;  \ell \ge 0; \\
\end{array}
\end{equation}

{\rm (iii)}~  up to topological isomorphism  in the category $\fr$ under the assumption that  $\A$ and $\B$  are nuclear and
up to isomorphism of linear spaces in the category of complete nuclear $DF$-spaces and in the category $\fr$ under the assumption that either $\A$ or  $\B$ nuclear, 
\begin{equation}
\begin{array}{ccccccccccc} \label{nice-HP-cohomology-AB}
{\mathcal H}{\mathcal P}_0(\A\ptp \B) \iso
 \A/[\A, \A] \ptp \B/[\B,\B] \;{\rm and} \;
{\mathcal H}{\mathcal P}_1(\A) \iso \{0\};\\
{\mathcal H}{\mathcal P}^0(\A\ptp \B) \iso (\A\ptp \B)^{tr}
 \;{\rm and} \;
{\mathcal H}{\mathcal P}^1(\A\ptp \B) \iso \{0\}.\\
\end{array}
\end{equation}
\end{corollary}

\begin{proof} Since $\A$ and $\B$ are unital,
${\mathcal H}^{bar}_n(\A \ptp \B) \iso \{0\}\;$ for  all  $n \ge 0$.
By Theorem \ref{A-X-trivial-Kunn-hom}, up to topological isomorphism,
$${\mathcal H}^{naive}_{0}(\A \ptp \B) \iso \A/[\A, \A] \ptp \B/[\B,\B]$$
and so is Hausdorff, and 
 ${\mathcal H}^{naive}_{n}(\A \ptp \B) \iso \{0 \}$ for  all  $n \ge 1$.
The result follows from \cite[Theorem 5.3]{Ly5}.
Note that, by definition, the  Hausdorff 
$~{\mathcal H}^{naive}_{0}(\A \ptp \B) \iso
(\A \ptp \B)/[\A\ptp \B, \A\ptp \B]$.
\end{proof}


\end{document}